\documentclass[a4paper,11pt]{article} 
\usepackage{amsmath, amssymb, esint, comment}
\usepackage{amscd}
\usepackage{shadow, epsfig, color}
\usepackage{graphicx}
\usepackage{url}
\usepackage{mathrsfs}
\usepackage[final]{changes_0}
\usepackage[normalem]{ulem}
\usepackage{xcolor,cancel}

\usepackage[hidelinks,
colorlinks=false,
linkcolor=blue,
filecolor=blue,
urlcolor=cyan,
citecolor=magenta
]{hyperref}

\newcommand{\N}{\mathbb N}

\newcommand{\Q}{\mathbb Q}
\newcommand{\R}{\mathbb R}
\newcommand{\C}{{\sf Ch}}

\newcommand{\e}{\varepsilon}
\newcommand{\1}{\mathbf 1}
\newcommand{\qed}{\ \hfill \fbox{} \bigskip}
\newcommand{\proof}{{\it Proof.} }
\newcommand{\dis}{\displaystyle}

\newcommand{\E}{\mathcal E}

\newcommand{\F}{\mathcal F}
\newcommand{\x}{\overline{x}}

\newcommand{\EN}{\overline{\N}}

\newcommand{\B}{\mathbb B}
\newcommand{\la}{\langle}
\newcommand{\ra}{\rangle}
\newcommand{\m}{\widetilde{m}}

\renewcommand{\tilde}{\widetilde}
\renewcommand{\comment}{\comment}

\numberwithin{equation}{section}

\newtheorem{thm}{Theorem}[section]
\newtheorem{defn}[thm]{Definition}
\newtheorem{lem}[thm]{Lemma}
\newtheorem{rem}[thm]{Remark}

\newtheorem{exa}[thm]{Example}
\newtheorem{prop}[thm]{Proposition}
\newtheorem{asmp}[thm]{Assumption}

\makeatletter
\newcommand{\subjclass}[2][2010]{%
  \let\@oldtitle\@title%
  \gdef\@title{\@oldtitle\footnotetext{#1 \emph{Mathematics Subject Classification.} #2}}%
}
\newcommand{\keywords}[1]{%
  \let\@@oldtitle\@title%
  \gdef\@title{\@@oldtitle\footnotetext{\emph{Key Words and Phrases.} #1.}}%
}
\makeatother

\renewcommand{\tilde}{\widetilde}

\title{Convergence of Brownian Motions on Metric Measure Spaces Under Riemannian Curvature--Dimension Conditions}
  \author{Kohei Suzuki \vspace{2mm} \\ {\it Scuola Normale Superiore
} \vspace{2mm}
   \\ {\it \small Piazza dei Cavalieri, 7, 56126, Pisa Italy}
            \\ {\small E-mail: kohei.suzuki@sns.it} \vspace{3mm} 
            }

\subjclass{Primary 60F17; Secondary 53C23.}
\keywords{Riemannian Curvature-Dimension Condition, Measured Gromov--Hausdorff Convergence, Brownian Motion, Weak Convergence}

\begin{document}
\date{}
\maketitle
\begin{abstract}
We show that the pointed measured Gromov convergence of the underlying spaces implies (or under some condition, is equivalent to) the weak convergence of Brownian motions under Riemannian Curvature-Dimension conditions. 
\deleted{This paper is an improved and jointed version of the previous two manuscripts \cite{S15} and \cite{S16}. The improvements extend our results to the case of $\sigma$-finite reference measures and to the case that initial distributions of Brownian motions are possible to be dirac measures.}
\end{abstract}
\tableofcontents

\section{Introduction}
\subsection{Motivation}
The aim of this paper is to characterize a probabilistic convergence of Brownian motions in terms of a geometric convergence of the underlying spaces. Our main results show that the pointed measured Gromov (pmG) convergence of the underlying spaces implies (or under some condition, is equivalent to) the weak convergence of Brownian motions under Riemannian Curvature-Dimension (RCD) conditions for the underlying spaces. 

Let us consider the following motivating example: let a sequence of Riemannian manifolds $\{M_n\}_{n\in \N}$
converges to a (possibly non-smooth) metric measure space in the Gromov--Hausdorff (GH) sense. Let $(B^n, \mathbb P^n)$ be a Brownian motions on each $M_n$. Noting that $(B^n, \mathbb P^n)$ can be determined only by the underlying geometric structure of the Riemannian manifolds $M_n$, an important question is whether 
\begin{enumerate}
\item[(Q)] a sequence of Brownian motions on Riemannian manifolds also converges weakly to the Brownian motion on the GH-limit space\deleted{, or not}. 
\end{enumerate}
This question however does not make sense without additional assumptions because there is a gap between the geometric and probabilistic convergences: the weak convergence of Brownian motions clearly involves the first-order differentiable structure of the underlying spaces although the GH convergence never sees any information of differentiable structures. Indeed,  we have examples whereby the limit process is no more a diffusion process (see (ii) and (iii) in Remark \ref{rem: relw}).

In this paper, adopting as an additional assumption {\it the uniform lower Ricci curvature bound of $M_n$}, we can answer (Q) affirmatively, which is an application of the main results in this paper.
To be more precise, we obtain the equivalence between these geometric/probabilistic convergences in the framework of metric measure spaces under {\it the synthetic lower Ricci curvature bound} (called RCD in this paper).

Let us explain the background issues in more detail. Generally, the GH-limit spaces of Riemannian manifolds with lower Ricci curvature bounds (called {\it Ricci limit spaces}) are so singular that they are not necessarily even topological manifolds and may have a dense singular set (see Example \ref{exa: dens}). However, they still have ``Riemannian-like" structures and similar properties to smooth Riemannian manifolds with lower Ricci curvature bounds, which have been investigated initially by Cheeger--Colding \cite{CC97, CC00a, CC00b}. 
The RCD condition
  is a proper generalization of the notion of lower Ricci curvature bounds to non-smooth spaces including Ricci limit spaces (see \cite{AGS14b, AGMR15, AMS16, G15, EKS15}). It is known that RCD spaces include various finite- and infinite-dimensional non-smooth spaces, not only Ricci limit spaces, but also  infinite-dimensional spaces such as Hilbert spaces with log-concave measures (related to various stochastic partial differential equations) (see further details in Section \ref{sec: exa}). 

By recent developments of analysis on metric measure spaces, we can construct Brownian motions on RCD spaces by using a certain quadratic form, what is called {\it Cheeger energy}. This is a generalization of Dirichlet energy on smooth manifolds and induces a quasi-regular strongly local conservative symmetric Dirichlet form (Ambrosio-Gigli-Savar\'e \cite{AGS14, AGS14b}, \added{Ambrosio--Gigli--Mondino--Rajala} \cite{AGMR15}), which is determined only by the underlying metric measure structure. 

One of the important problems for Brownian motions on these non-smooth spaces is to characterize the weak convergence of Brownian motions in terms of some geometrical convergence of the underlying spaces, which we call the {\it stability} of Brownian motions. The significance of the stability can be explained from several different perspectives. From the standpoint of limit theorems of stochastic processes, the stability is interpreted as a geometric characterization of invariance principles for Brownian motions in the sense that  Brownian motions on limit spaces are approximated by Brownian motions on converging spaces.  
From the viewpoint of ``well-definedness", the stability also enables us to verfiy that  Brownian motions in limit spaces are ``well-defined" in the sense that Brownian motions intrinsically defined by Cheeger energies on limit spaces coincide with limit processes of Brownian motions on approximating spaces. 
From the perspective that  Brownian motions are  considered as ``a map" assigning laws of diffusions (i.e., probability measures on path spaces) to each metric measure space, the stability reveals the interesting fact that this map is continuous with respect to the corresponding topologies (e.g., GH-topology of metric measure spaces/weak topology of probability measures on path spaces), which is one ideal aspect of Brownian motions but has not been focused on so much until now. 

The main contribution of this paper is to prove the stability of Brownian motions in the general framework of RCD spaces, whereby various singular/infinite-dimensional spaces are included. Moreover, we show several equivalences of the weak convergence of Brownian motions and the pmG convergence of the underlying spaces. \deleted{We remark that the pmG convergence of the underlying spaces is obviously not enough to imply the weak convergence of Brownian motions without any additional condition, which is, in this paper, the RCD condition. This is because the weak convergence of Brownian motions involves first-order differential structures since the Cheeger energy includes information of the gradient in its definition, but the pmG convergence never sees any kind of differential structures. By aid of the RCD condition, which plays a role as a control of a second-order differential structure (Ricci curvature), we obtain the stability of Brownian motions.}For references to other investigations regarding the stability problem, see the historical remarks (Section \ref{sec: history} below).

\subsection{Main Results}
In this paper, we always consider {\it pointed metric measure (p.m.m.)\ spaces} $\mathcal X=(X, d, m, \x)$ whereby 
\begin{align} \label{asmp: basic}
&\text{$(X,d)$ is a complete separable geodesic metric space with nonnegative and nonzero} \notag
\\
&\ \text{Borel measure $m$ which is finite on all bounded sets, and $\x$ is a fixed point in ${\rm supp}[m]$.} 
\end{align}
\deleted{Our main results consist of two parts, one is for RCD$(K,\infty)$ spaces, and the other is for RCD$^*(K,N)$ spaces. The latter condition is stronger than the former one.
We first state the results for  RCD$(K,\infty)$ spaces.}
For\deleted{the} main theorems, we assume the following condition:
\begin{asmp} \label{asmp: 2} \normalfont Let $K \in \R$ and $\EN:=\N \cup \{\infty\}.$
Let $\{\mathcal X_n \}_{n \in \EN}=\{(X_n, d_n, m_n, \x_n)\}_{n \in \EN}$ be a sequence of p.m.m.\ spaces satisfying \eqref{asmp: basic} and RCD$(K,\infty)$ condition.
\end{asmp}
The notion of ${\rm CD}(K,\infty)$  spaces was introduced by Sturm \cite{Sturm06} and Lott--Villani \cite{LV09}, and the notion of $\mathrm{RCD}(K,\infty)$ spaces was introduced by Ambrosio--Gigli--Savar\'e \cite{AGS14b} and Ambrosio--Gigli--Mondino--Rajala \cite{AGMR15}. The CD$(K,\infty)$ condition is a generalization of Ricci curvature bounded from below by $K$ to metric measure spaces in terms of the $K$-convexity of the entropy on the Wasserstein spaces. Furthermore RCD$(K,\infty)$ condition means the CD$(K,\infty)$ and that the Cheeger energy is {\it quadratic}. We will explain the precise definition in Subsection \ref{subsec: RCD}. RCD spaces admit the GH limit spaces of Riemannian manifolds with lower Ricci curvature bounds, and also admit Alexandrov spaces (metric spaces satisfying a generalized notion of {\rm ``sectional curvature}$\ge K"$) (Petrunin \cite{Pet11} and Zhang--Zhu \cite{ZZ10})\added{, cone spaces and warped product spaces (Ketterer \text{$\cite{Ket14, Ket14a}$}), and quotient spaces (Galaz-Garc\'ia--Kell--Mondino--Sosa \cite{GKMS17})}. Moreover, not only finite-dimensional spaces, but also  several infinite-dimensional spaces related to stochastic partial differential equations are included such as Hilbert spaces with log-concave measures (Ambrosio--Savar\'e--Zambotti \cite{ASZ09}).

Under Assumption \ref{asmp: 2}, we can always take constants $c_1, c_2 >0$ 
satisfying the following volume growth estimate (see \cite[Theorem 4.24]{Sturm06})
\begin{align} \label{asmp: basic1}
m_n(B_r(\x_n)) \le c_1e^{c_2r^2}, \quad \forall r>0.
\end{align}
Here we mean $B_r(\x_n):=\{x \in X_n: d(x,\x_n)<r\}$. Taking ${\sf C}>c_2$, we set {\it a weighted measure $\m_n$} as follows:
\begin{align} \label{constant: vol}
z_n:=\int_{X_n} e^{-{\sf C}d_n^2(x,\x_n)}dm_n(x), \quad \text{and} \quad \m_n:=
\begin{cases} \dis
\frac{1}{z_n}e^{-{\sf C}d_n^2(\cdot, \x_n)}m_n, &\text{if $m_n(X_n)=\infty$,}
\\ \dis
\frac{1}{m_n(X_n)}m_n, &\text{if $m_n(X_n)<\infty$}.
\end{cases}
\end{align}
Under Assumption \ref{asmp: 2},  the Cheeger energy ${\sf Ch}_n$ on $\mathcal X_n=(X_n,d_n,m_n, \x_n)$ (see Subsection \ref{subsub: Ch}) induces a quasi-regular conservative symmetric strongly local Dirichlet form, and there exists a conservative \replaced{symmetric Markov}{diffusion} process $(\{\mathbb P_n^{x}\}_{x \in X_n}, \{B_t^n\}_{t \ge 0})$ \deleted{$(\{\mathbb P_n^{x}\}_{x \in X_n}, \{B_t^n\}_{t \ge 0})$} on $\mathcal X_n$ (See Section \ref{subse: BM1}).
We call $(\{\mathbb P_n^{x}\}_{x \in X_n}, \{B_t^n\}_{t \ge 0})$ {\it Brownian motion} on $\mathcal X_n$.
 
The \added{following} main theorem states that the weak convergence of the Brownian motions can be characterized by the pmG convergence of the underlying spaces under Assumption \ref{asmp: 2} (we will give the definition of the pmG convergence in Subsection \ref{subsec: D}). 
\begin{thm} \label{thm: mthm2}
 Suppose that Assumption \ref{asmp: 2} holds. Then the following (i) and (ii) are equivalent:
 \begin{itemize}
 \item[{\bf (i)}] {\rm \bf \small (pmG Convergence of the Underlying Spaces)} \vspace{2mm} \\
The p.m.m.\ spaces $\{\mathcal X_n\}_{n \in \N}$ converge to $\mathcal X_\infty=(X_\infty, d_\infty, m_\infty, \x_\infty)$ in the pmG sense.
 \item[{\bf (ii)}] {\rm \bf \small (Weak Convergence of the Laws of Brownian Motions)} \vspace{2mm}\\  
There exist a complete separable metric space $(X,d)$ and isometric embeddings $\iota_n: X_n \to X\ (n \in \EN)$ so that $\iota_n(\x_n) \to \iota_\infty(\x_\infty)$, and 
\begin{align*} 
(\iota_n(B^n), \mathbb P_n^{\m_n}) \ {\to} \ (\iota_\infty(B^\infty), \mathbb P_\infty^{\m_\infty}) \quad \text{weakly} \quad \text{in $\mathcal P(C([0,\infty); X))$}.
\end{align*}
\end{itemize}
\end{thm}
Here $(\iota_n(B^n), \mathbb P_n^{\m_n})$ means the law of the embedded Brownian motion $\iota_n(B^n)$ with the initial distribution $\m_n$ and $\mathcal P(C([0,\infty); X))$ denotes the set of all Borel probability measures on the continuous path space $C([0,\infty); X)$.
\begin{rem} \label{rem: AC} \normalfont 
Several remarks for Theorem \ref{thm: mthm2} are given below.				
\begin{enumerate}
\item[(i)] 
The ${\rm RCD}(K,\infty)$ condition is stable under the pmG convergence (see \cite[Theorem 7.2]{GMS13}), and therefore the limit space $\mathcal X_\infty$ also satisfies the ${\rm RCD}(K,\infty)$ condition so that the Brownian motion can be defined also on the limit space $\mathcal X_\infty$.
\item[(ii)] The pmG convergence is weaker than the measured Gromov-Hausdorff convergence. See \cite[Theorem 3.30]{GMS13}.
\end{enumerate}
\end{rem}

In statement (ii) in Theorem \ref{thm: mthm2}, the initial distribution is absolutely continuous with respect to the reference measure $m_n$. It is natural in the next step to ask how the case of the dirac measure $\delta_{\x_n}$ is, which means the Brownian motions start at the point $\x_n$. We introduce several conditions below:
 \begin{enumerate}
 \item[{\bf (A)}] For any $n \in \N$, $m_n(X_n)=1$.
 \item[{\bf (B)}] For any $r>0$ \added{and any $t>0$}, 
 $$\sup_{n \in \N}\|p_n(t,\x_n, \cdot)\|_{\infty, B_r(\x_n)}<\infty,$$
 whereby $p_n(t,x,y)$ is the density of the transition probability $p_n(t,x,dy)$ of $\{P_t^n\}_{t \ge 0}$ with respect to the reference measure $m_n$, and $\|\cdot\|_{\infty, B_r(\x_n)}$ means the essential supremum on the ball $B_r(\x_n)$.
 \end{enumerate}
Now we state the second main result. 
\begin{thm} \label{thm: mthm2-0}
 Suppose that Assumption \ref{asmp: 2} holds.
 If, moreover, either  {\bf (A)}, or {\bf (B)} holds,
 then (i) (thus also (ii))  in Theorem \ref{thm: mthm2} implies the following (iii)$_{\ge \e}$: for any $\e>0$, 
\begin{itemize}
 \item[{\bf (iii)}$_{ \ge \e }$] {\rm \bf \small (Weak Convergence of the Laws of Brownian Motions \replaced{S}{s}tarting at \replaced{P}{p}oints in a Time Interval $[\e,\infty)$)} \vspace{2mm}\\  
 There exist a complete separable metric space $(X,d)$ and isometric embeddings $\iota_n: X_n \to X\ (n \in \EN)$ so that $\iota_n(\x_n) \to \iota_\infty(\x_\infty)$ and 
 \begin{align} \label{main: WEAKCON2}
(\iota_n(B^n), \mathbb P_n^{\x_n}) \ {\to} \ (\iota_\infty(B^\infty), \mathbb P_\infty^{\x_\infty}) \quad \text{weakly} \quad \text{in $\mathcal P(C([\e,\infty); X))$}.
\end{align}
\end{itemize}
\end{thm}
\begin{rem} \label{rem: AC} \normalfont 
Several remarks for Theorem \ref{thm: mthm2-0} are given below.				
\begin{enumerate}
\item[(i)] Condition {\bf (B)} is satisfied for any RCD$^*(K,N)$ spaces \replaced{according to the Gaussian heat kernel estimate by Jiang--Li--Zhang \cite{JLZ15}.}{because of local Gaussian heat kernel estimates, which follow from the local volume doubling property and the local Poincar\'e inequality according to Sturm \cite{Sturm96} (see also Jiang--Li--Zhang \cite{JLZ15}).}
\item[(ii)]
If the following uniform ultra-contractivity of the heat semigroup $\{H_t\}_{t \ge 0}$ holds, then  condition {\bf (B)} holds (see {\cite[Proposition 6.4]{AGS14b}}): there exists a $p>1$ so that, with some positive constant  $C(t,K)$ dependent only on $t$ and $K$, we have
$$\|H_tf\|_p \le C(t,K)\|f\|_1, \quad \forall f \in L^1(X,m),\ \forall t>0.$$
We have examples satisfying the ultra-contractivity which is a RCD$(K,\infty)$ space but not a RCD$^*(K,N)$ space for any $1<N<\infty$. 
Let $\mathcal X_\alpha=(\R, |\cdot-\cdot |, C_\alpha\exp\{-|\cdot|^\alpha\}dx)$ whereby $\alpha \in \{2,4,6,...\}$ is an even number and $C_\alpha$ is the normalizing constant. For any $\alpha>2$, it is known that $\mathcal X_\alpha$ satisfies the ultra-contractivity of the heat semigroup (Kavian--Kerkyacharian--Roynette \cite{KKR93}) and satisfies the RCD$(0,\infty)$ condition, but not RCD$^*(K,N)$ for any finite $1<N<\infty$.
\end{enumerate}
\end{rem}
The notion of $\mathrm{RCD}^*(K,N)$ condition is a generalization of Ricci curvature bounded from below by $K$ and dimension bounded above by $N$ to metric measure spaces, which is stronger than the RCD$(K,\infty)$ condition (see \cite{G15, EKS15, AMS16}).

Next we consider the converse implication that the weak convergence of Brownian motions induces the pmG convergence of the underlying spaces. 
\deleted{Let $p_n(t,x,y)$ be the heat kernel for the p.m.m.\ space ${\mathcal X}_n=(X_n, d_n, m_n, \x_n)$}\added{Define $p_n(t, x,x)=\|p_n(t/2, x, \cdot)p_n(t/2, x, \cdot)\|^2_2$, which can be defined for every $x \in X$}. Let us consider the following condition: there exists $t_* >0$ and a constant $M$ so that
\begin{align} \label{eq: bound of HK} 
\sup_{n \in \N}p_n(t_*,\x_n, \x_n)<M<\infty. 
\end{align}
Note that, since $p_n(t,x,x)$ is non-increasing function in $t$, if we find the time $t_*$ satisfying \eqref{eq: bound of HK}, then for any $t>t_*$, the estimate \eqref{eq: bound of HK} holds. 
\replaced{We also note that,}{For instance,} \replaced{for fixed $1<N<\infty$ and $K \in \R$, if $X_n$ satisfies the RCD$^*(K,N)$ for all $n \in \N$ and $\inf_n m_n(B_{1}(\x_n))>0$ (this holds under (i), (ii), (iii)$_{\ge \e}$, or (iii)$_{\ge 0}$ since $\{\iota_n(\x_n)\}_{n \in \N}$ is bounded in a common ambient space $X$ according to $\iota(\x_n) \to \iota_\infty(\x_\infty)$), then \eqref{eq: bound of HK} is satisfied for some constant $M$ and $t_*$ because of the local Gaussian heat kernel estimate by Jiang-Li-Zhang \cite{JLZ15}.}{all RCD$^*(K,N)$ spaces with $1<N<\infty$ satisfies \eqref{eq: bound of HK}  because of the local Gaussian heat kernel estimate.}
Let ${\rm diam}(X_n)$ denote the diameter of $X_n$: ${\rm diam}(X_n):=\sup_{x,y \in X_n}d_n(x,y)$. 
 We now state the following theorem:
\begin{thm} \label{thm: mthm2-1}
Suppose that Assumption \ref{asmp: 2} and{\deleted{the} condition \eqref{eq: bound of HK}} hold. {If, m}oreover, either $K>0$, or $\sup_{n \in \N}{\rm diam}(X_n)<D$ holds \added{for some $0<D<\infty$},
then (iii)$_{ \ge \e}$ for any $\e>0$ in Theorem \ref{thm: mthm2-0} implies  (i) and (ii) in Theorem \ref{thm: mthm2} $($therefore all statements (i), (ii) and (iii)$_{\ge \e}$ for any $\e>0$ are equivalent$)$. 
\end{thm} 
Finally, we give the following statement, in which all statements (i), (ii), (iii)$_{\ge \e}$ for any $\e>0$, and (iii)$_{\ge 0}$ are equivalent under the RCD$^*(K,N)$ condition with a uniform diameter bound.
\begin{thm} \label{thm: mthm1} Let $K \in \R$, $1<N<\infty$ and $0<D<\infty$.
Suppose that a sequence of p.m.m.\ spaces $\{\mathcal X_n\}_{n \in \N}$ satisfies \eqref{asmp: basic}, RCD$^*(K,N)$ and $\sup_{n \in \N}{\rm diam}(X_n)<D$. Then all\deleted{the} four statements of (i), (ii) in Theorem \ref{thm: mthm2}, (iii)$_{ \ge \e}$ for any $\e>0$ in Theorem \ref{thm: mthm2-0} and the following $(iii)_{ \ge 0}$ are  equivalent:
 \begin{itemize}
 \item[{\bf (iii)}$_{ \ge 0}$] {\rm \bf \small (Weak Convergence of the Laws of Brownian Motions \replaced{S}{s}tarting at \replaced{P}{p}oints  in a Time Interval $[0,\infty)$)} \vspace{2mm}\\  
There exist a compact metric space $(X,d)$ and isometric embeddings $\iota_n: X_n \to X\ (n \in \EN)$ so that
\begin{align*} 
(\iota_n(B^n), \mathbb P_n^{\x_n}) \ {\to} \ (\iota_\infty(B^\infty), \mathbb P_\infty^{\x_\infty}) \quad \text{weakly} \quad \text{in $\mathcal P(C([0,\infty); X))$}.
\end{align*}
\end{itemize}
\end{thm}
\begin{rem} \label{rem: AC1} \normalfont 
We give several remarks for Theorem \ref{thm: mthm1}.
\begin{enumerate}
\item[(i)] The ${\rm RCD}^*(K,N)$ condition is stable under the pmG convergence (see \cite{EKS15}), and therefore the limit space $\mathcal X_\infty$ also satisfies the ${\rm RCD}^*(K,N)$ condition so that the Brownian motion can be defined at every starting point also on the limit space $\mathcal X_\infty$.
\item[(ii)] The pmG convergence (see Definition \ref{prop: Dconv}) is equivalent to the pointed measured Gromov-Hausdorff convergence under \replaced{the assumption in Theorem \ref{thm: mthm1}}{Assumption \ref{asmp: 1}}  (see \cite[Theorem 3.33]{GMS13}).
\end{enumerate}
\end{rem}

\subsection{Historical Remarks} \label{sec: history}
\begin{rem} \label{rem: relw}\normalfont 
Several historical remarks are given below.
\begin{description}
\item[(i)]  In Ambrosio--Savar\'e--Zambotti \cite[Theorem 1.5]{ASZ09},\deleted{under the weak convergence of the underlying reference measures,} they investigated the weak convergence of Brownian motions on \added{a fixed} Hilbert space\deleted{s} (as an ambient space) with \added{varying} log-concave measures and norms, which is a specific case of RCD$(0,\infty)$ spaces. \added{Their metrics $d_n$ are not necessarily isometric to the metric $d$ in the ambient space, but each $d_n$ is equivalent to $d$.} \added{In the case that each $d_n$ is isometric to $d$,} our results (Theorem \ref{thm: mthm2} and \ref{thm: mthm2-0}) \replaced{can be seen as a generalization of}{generalize} their result \cite[Theorem 1.5]{ASZ09} \replaced{to}{for} general RCD$(K,\infty)$ spaces.

\item[(ii)] In Ogura \cite{O01}, under the condition of uniform upper bounds for heat kernels (not necessarily lower bound of Ricci curvatures) and the Kasue-Kumura (KK) spectral convergence, he studied the weak convergence of the laws of {\it time-discretized} Brownian motions on weighted compact Riemannian manifolds. The KK spectral convergence roughly means a uniform convergence of heat kernels and stronger than the mGH convergence. In his case, the Ricci curvature is not necessarily bounded from below and the limit process may be a jump process (\cite[4.6]{O01}). 
The time-discretization is one possible approach for a convergence of stochastic processes on varying spaces, while we adopt in this paper a different approach, i.e., embedding into one common metric space $X$. 
\item[(iii)] If we do not assume RCD conditions for a sequence of the underlying metric measure spaces, then limit processes are not necessarily diffusions. In Ogura--Tomisaki--Tuchiya \cite{OTT02}, they considered a sequence of Euclidean spaces $(\R^d, \|\cdot\|_2)$ with certain underlying measures $\mu_n$ whereby $\{(\R^d, \|\cdot\|_2, \mu_n)\}_{n \in \N}$ \replaced{do}{does} not necessarily satisfy RCD conditions. They showed that diffusion processes on $\R^d$ associated with the corresponding local Dirichlet forms converge to jump processes (or generally jump-diffusion processes) corresponding to certain non-local Dirichlet forms. 

\item[(iv)]
In Freidlin--Wentzell \cite{FW93, FW04} and Albeverio-Kusuoka \cite{AK12} (see also references therein), diffusion processes associated with SDEs on thin tubes in $\R^d$ were studied. When thin tubes shrink to  a spider graph, diffusion processes converge weakly to  a one-dimensional diffusion on this spider graph.
Their setting does not satisfy the $\mathrm{RCD}$ condition since spider graphs branch at points of conjunctions but {RCD spaces are essentially non-branching (see \cite[{Theorem 1.1}]{RS14}).}

\item[(v)]
In Athreya--L\"ohr--Winter \cite{ALW}, the weak convergence of certain Markov processes on tree-like spaces was studied. When tree-like spaces converge in {\it Gromov-vague} sense, the corresponding processes also converge weakly. Their tree-like spaces admit 0-hyporbolic spaces, which are not necessarily included in RCD spaces. 

\item[(vi)] In Suzuki \cite{S15}, the author investigated the weak convergence of continuous stochastic processes on metric spaces converging in the Lipschitz distance. The Lipschitz convergence is stronger than the measured Gromov convergence (see \cite[Section 3.C]{Gro98}).
\end{description}
\end{rem}

Finally we list related studies not mentioned in Remark \ref{rem: relw}.
\added{In Suzuki \cite{S17}, the author studied the weak convergence of non-symmetric diffusion processes on RCD spaces as a next step of the current paper.}
In Li \cite{L16, L17}, she studied a convergence of random ODE/SDE on manifolds. In Stroock--Varadhan \cite{SV79}, Stroock--Zheng \cite{SZ97} and Burdzy--Chen \cite{BC08}, approximations of diffusion processes on $\mathbb R^d$ by discrete Markov chains on $(1/n)\mathbb Z^d$ were investigated. In Bass--Kumagai--Uemura \cite{BKU10} and Chen--Kim--Kumagai \cite{CKK13}, they studied approximations of jump processes on proper metric spaces by Markov chains on discrete graphs. Approximations of Markov processes on ultra-metric spaces were explored in Suzuki \cite{S14}. In Pinsky \cite{Pin76}, he studied approximations of Brownian motions on Riemannian manifolds by random walks, while the case of sub-Riemannian manifolds was investigated by Gordina and Laetsch \cite{GL14}. 
In Croydon--Hambly--Kumagai \cite{CHK16}, in which it was assumed that a sequence of resistance forms converges with respect to the GH-vague topology and satisfies a uniform volume doubling condition, they showed the weak convergence of corresponding Brownian motions and local times.
There are many studies about scaling limits of random processes on random environments (see, e.g., Kumagai \cite{K14} and references therein).

\subsection{Organization of the Paper}
The paper is structured as follows: First, the notation is fixed and preliminary facts are recalled in Section \ref{sec: Pre} (no new results are included)\deleted{.Each of the subsections define important aspects of the research in this paper}, namely: basic notations and basic definitions (Subsection \ref{subsec: Pre}); $L^2$-Wasserstein distance (Subsection \ref{subsec: L2}); pmG convergence (Subsection \ref{subsec: D}); RCD$(K,\infty)$ and RCD$^*(K,N)$ spaces (Subsection \ref{subsec: RCD}); $L^2$-convergence of the heat semigroup (Subsection \ref{sec: MC}).
In Section \ref{sec: BM}, we state several properties about Brownian motions on RCD spaces. 
In Section \ref{sec: exa} we present examples in which Assumptions \ref{asmp: 2} and the assumption in Theorem \ref{thm: mthm1} are satisfied. These examples consist of weighted Riemannian manifolds and \replaced{their}{its} pmG limit spaces, Alexandrov spaces, and Hlibert spaces with log-concave probability measures.  
In Section \ref{sec: proof2}, we give the proof of Theorem \ref{thm: mthm2} . In Section \ref{subsec: I-III}, we show the proof of Theorem \ref{thm: mthm2-0}. In Section \ref{sec: proof}, we prove Theorem \ref{thm: mthm2-1}. Finally, in Section \ref{sec: proof1}, we prove Theorem  \ref{thm: mthm1}.

\section{Notation \& Preliminary Results} \label{sec: Pre}
\subsection{Notation} \label{subsec: Pre}
Let $\N=\{0,1,2,...\}$ and $\overline{\N}:=\N \cup \{\infty\}$ denote the set of natural numbers and the set of natural numbers with $\{\infty\}$ respectively.  
For a complete separable metric space $(X,d)$, we denote by $B_r(x)=\{y \in X: d(x,y)<r\}$ the open ball centered at $x \in X$ with radius $r>0$.
By using $\mathscr B(X)$, we mean the family of all Borel sets in $(X,d)$; and by $\mathcal B_b(X)$, the set of real-valued bounded Borel-measurable functions on $X$.
Let $C(X)$ be the set of real-valued continuous functions on $X$, while $C_b(X)$, \added{$C_{\infty}(X)$,}  $C_0(X)$ and $C_{bs}(X)$ denote the subsets of $C(X)$ consisting of bounded functions, \added{functions vanishing at infinity,} functions with compact support, and bounded functions with bounded support, respectively. \added{Let ${\rm Lip}(X)$ and ${\rm Lip}_b(X)$ denote the set of Lipschitz continuous functions, and the set of bounded Lipschitz continuous functions, respectively. For $f \in {\rm Lip}(X)$, we denote by ${\rm Lip}_X(f)$ the global Lipschitz constant of $f$.}
The set $\mathcal P(X)$ denotes all Borel probability measures on $X$.
The set of continuous functions on $[0,\infty)$ valued in $X$ is denoted by $C([0,\infty), X)$. 

A continuous curve $\gamma: [a,b] \to X$ is {\it connecting $x$ and $y$} if $\gamma_a=x$ and $\gamma_b=y$.
A continuous curve $\gamma: [a,b] \to X$ is  {\it a minimal geodesic} if 
$$d(\gamma_t, \gamma_s)=\frac{|s-t|}{|b-a|}d(\gamma_a, \gamma_b) \quad a \le t \le s \le b.$$
In particular, if $\frac{d(\gamma_a, \gamma_b)}{|b-a|}$ can be replaced by $1$, we say that $\gamma$ is {\it  unit-speed}. \added{A metric space $X$ is called {\it geodesic} if for any two points $x,y\in X$, there exists a minimal geodesic $\{\gamma_t\}_{t \in [0,1]}$ connecting $x$ and $y$.}

Let $\mathrm{supp}[m]=\{x \in X: m(B_r(x))>0, \ \forall r>0\}$ denote the support of $m$.
Let $(Y,d_Y)$ be a complete separable metric space. For a Borel measurable map $f: X \to Y$, let $f_\#m$ denote the push-forward measure on $Y$:
${f}_\#m(B) = m(f^{-1}(B))$ for any Borel set $B \in \mathscr B(Y).$

\subsection{$L^p$-Wasserstein Space} \label{subsec: L2}
Let $(X_i,d_i)$ $(i=1,2)$ be complete separable metric spaces and $1\le p<\infty$.
For $\mu_i \in \mathcal P(X_i)$, a probability measure $q \in \mathcal P(X_1 \times X_2)$ is called {\it a coupling of $\mu_1$ and $\mu_2$} if 
${\pi_1}_\#q = \mu_1$ and ${\pi_2}_\#q=\mu_2,$
whereby $\pi_i$ $(i=1,2)$ is the projection $\pi_i: X_1 \times X_2 \to X_i$ as $(x_1,x_2) \mapsto x_i$. We denote by $\Pi(\mu, \nu)$ the set of all coupling of $\mu$ and $\nu$.

Let $(X,d)$ be a complete separable metric space. Let $\mathcal P_p(X)$ be the subset of $\mathcal P(X)$ consisting of all Borel probability measures $\mu$ on $X$ with finite $p$-th moment:
$$\int_Xd^p(x,\x)d\mu(x)<\infty \quad \text{for some (and thus any) $\x \in X.$}$$
We equip $\mathcal P_p(X)$ with the\deleted{quadratic} transportation distance $W_p$, called {\it $L^p$-Wasserstein distance}, defined as follows:
\begin{align} \label{eq: Was}
W_p(\mu,\nu)=\Bigl(\inf_{q \in \Pi(\mu,\nu)}\int_{X\times X}d^p(x,y)dq(x,y)\Bigr)^{1/p}.
\end{align}
A coupling $q \in \Pi(\mu,\nu)$ is called {\it an optimal coupling} if $q$ attains the infimum in the equality \eqref{eq: Was}. It is known that, for any $\mu,\nu$, there always exists an optimal coupling $q$ of $\mu$ and $\nu$ (e.g., \cite[\S 4]{V09}).
It is also known that $(\mathcal P_p(X), W_p)$ is a complete separable \added{geodesic} metric space for $1<p<\infty$ \added{if $(X,d)$ is a complete separable geodesic metric space} (e.g., \cite[Theorem 6.18]{V09}).

\subsection{Pointed Measured Gromov Convergence} \label{subsec: D}
In this subsection, we recall the definition of pmG convergence introduced in Gigli-Mondino-Savar\'e \cite{GMS13}.
\begin{defn}[\cite{GMS13}]\normalfont ({\bf pmG Convergence}) \label{prop: Dconv} 
\\
A sequence of p.m.m.\ spaces $\{\mathcal X_n=(X_n, d_n, m_n, \x_n)\}_{n \in \N}$ satisfying \eqref{asmp: basic} is {\it convergent to $\mathcal X_\infty=(X_\infty, d_\infty, m_\infty, \x_\infty)$ in the pointed measured Gromov (pmG) sense} if there exist a complete separable metric space $(X,d)$ and isometric embeddings $\iota_n: X_n \to X\  (n \in \EN:=\N \cup \{\infty\})$
satisfying
\begin{align}\label{eq: VGC}
\iota_n(\x_n) \to \iota_\infty(\x_\infty) \in {\rm supp}[m_\infty], \quad \text{and}\quad \int_{X}f \ d({\iota_n}_\#m_n) \to \int_{X}f \ d({\iota_\infty}_\#m_\infty),
\end{align}
for any bounded continuous function $f:X \to \R$ with bounded support. 
\end{defn}
\begin{rem} \normalfont \label{rem: GHMG}
We would like to remark on the pmG convergence in Definition \ref{prop: Dconv}.
\begin{enumerate}
\item[(i)] In general, the pmG convergence is strictly weaker than the pointed measured Gromov-Hausdorff (pmGH) convergence (\cite[{Theorem} 3.30, Example 3.31]{GMS13}). However, if ${\rm supp}[m_\infty]=X_\infty$ and $\{\mathcal X_n\}_{n \in \N}$ satisfies a uniform doubling condition, then the two notions of pmG and pmGH coincide (see \cite[{Theorem 3.33}]{GMS13}).
\item[(ii)] The pmG convergence is metrizable by the distance $p\mathbb G_W$ on the collection $\mathbb X$ of all isomorphism classes of p.m.m.\ spaces (see \cite[Definition 3.13]{GMS13}). The space $(\mathbb X, p\mathbb G_W)$, moreover,  \replaced{is}{becomes} a complete and separable metric space (see \cite[Theorem 3.17]{GMS13}). 
\end{enumerate}
\end{rem}

\subsection{RCD Spaces} \label{subsec: RCD}
In this subsection, we recall the definition of the $\mathrm{RCD}(K,\infty)$ condition. 
We also recall several properties satisfied on $\mathrm{RCD}(K,\infty)$ spaces. See \cite{A18} for more comprehensive accounts of this field. 
\subsubsection{Relative Entropy}
In this subsection, we recall the definition of the relative entropy functional ${\rm Ent}_m: \mathcal P_2(X) \to \overline{\R}:=\R \cup \{+\infty\}$:
\begin{align*}
{\rm Ent}_m(\mu)=
\begin{cases} \dis
\int_X\frac{d\mu}{dm}\log(\frac{d\mu}{dm})dm, &\text{if $\mu \ll m$},
\vspace{3mm} \\ 
+\infty, &\text{otherwise}.
\end{cases}
\end{align*} 
Here $d\mu/dm$ denotes the Radon--Nikodym derivative. Let us write $D({\rm Ent}_m):=\{\mu \in \mathcal P_2(X): {\rm Ent}_m(\mu)<\infty\}$.
Although $m$ might not be a probability measure, the entropy ${\rm Ent}_m$ is well-defined and lower-semicontinuous thanks to condition \eqref{asmp: basic1}. Indeed, by recalling \eqref{constant: vol}:
\begin{align*}
z:=\int_X e^{-{\sf C}d^2(x,\x)}dm(x), \quad \text{so that} \quad \tilde{m}:=\frac{1}{z}e^{-{\sf C}d^2(\cdot, \x)}m,
\end{align*}
we can check that, for any $\rho m=\mu \in D({\rm Ent}_m)$ with $\rho=\frac{d\mu}{dm}$, it holds that $\mu=z\rho e^{{\sf C}d^2(\cdot,\x)}\tilde{m}.$ Therefore we obtain
\begin{align*} 
{\rm Ent}_m(\mu)={\rm Ent}_{\tilde{m}}(\mu)-{\sf C} \int_X d^2(\cdot, \x)d\mu - \log z,
\end{align*}
which implies that ${\rm Ent}_m$ is well-defined and lower-semicontinuous with respect to $W_2$-topology. 

\subsubsection{Cheeger Energy} \label{subsub: Ch} 
In this subsection, we recall the Cheeger energy ${\sf Ch}$ on $(X,d,m, \x)$. 
For $f \in \mathrm{Lip}(X)$, the local Lipschitz constant $|\nabla f|: X \to \R$ is defined as follows:
\begin{align*}|\nabla f|(x)&=
\begin{cases} \dis
\limsup_{y \to x}\frac{|f(y)-f(x)|}{d(y,x)}, & \text{if $x$ is not isolated},
\\ \dis
0, &\text{otherwise}.
\end{cases}
\end{align*}
Then we now recall the definition of {\it Cheeger energy}: (see \cite{G15, AGMR15, AGS14}) 
\begin{align*}
{\sf Ch}(f)&=\frac{1}{2}\inf\Bigl\{ \liminf_{n \to \infty} \int |\nabla f_n|^2 dm: f_n \in \mathrm{Lip}(X)\cap L^2(X,m), \int_X|f_n -f|^2 dm \to 0 \Bigr\}
\\
W^{1,2}(X,d,m)&=\{f \in L^2(X{,} m): {\sf Ch}(f)<\infty\}.
\end{align*}
If ${\sf Ch}(f) <\infty$, then the Cheeger energy can be written as an integral form by {\it minimal weak upper gradient $|\nabla f|_w$} (see \cite{AGS14, AGMR15}):
\begin{align*} 
{\sf Ch}(f)=\frac{1}{2}\int_X|\nabla f|^2_w dm, \quad \forall f \in W^{1,2}(X,d,m).
\end{align*}

\subsubsection{$\mathrm{RCD}(K,\infty)$ Spaces} 
In this subsection, we recall the $\mathrm{CD}(K,\infty)$/$\mathrm{RCD}(K,\infty)$ condition.
\begin{defn} \normalfont \label{defn: RCDCD}
The $\mathrm{CD}(K,\infty)$/$\mathrm{RCD}(K,\infty)$ conditions are defined as follows:
\begin{enumerate}  
\item[(i)] ($\mathrm{CD}(K,\infty)$) [\cite{Sturm06}, \cite{LV09}] 
\\
We say that $(X,d,m)$ satisfies {\it the curvature-dimension condition} $\mathrm{CD}(K,\infty)$ for $K \in \R$ if,
for each $\mu_0, \mu_1 \in D({\rm Ent}_m)$, there exists a $W_2$-geodesic $\{\mu_t\}_{t \in [0,1]} \subset D({\rm Ent}_m)$ connecting $\mu_0$ and $\mu_1$ so that 
\begin{align} \label{Ent}
{\rm Ent}_m({\mu_t}) \le (1-t){\rm Ent}_m(\mu_0) + t {\rm Ent}_m(\mu_1)-\frac{K}{2}t(1-t)W_2^2(\mu_0,\mu_1).
\end{align}
\item[(ii)] ($\mathrm{RCD}(K,\infty)$) [\cite[Remark 4.20]{G15}, \cite[Theorem 5.1]{AGS14b}, \cite[Theorem 6.1]{AGMR15}]
\\
We say that $(X,d,m)$ satisfies {\it the Riemannian curvature-dimension condition} $\mathrm{RCD}(K,\infty)$ if the following two conditions hold:
\begin{enumerate}
	\item[(ii-a)] $\mathrm{CD}(K,\infty)$
	\item[(ii-b)] the infinitesimal Hilbertianity, that is, the Cheeger energy {\sf Ch} is a quadratic form: 
	\begin{align*}
	2{\sf Ch}(u)+2{\sf Ch}(v)={\sf Ch}(u+v) + {\sf Ch}(u-v), 
	\end{align*}
	for any $u,v \in W^{1,2}(X,d,m).$
\end{enumerate}
\end{enumerate}
\end{defn}

It is known that CD$(K,\infty)$/RCD$(K,\infty)$ conditions are stable under the pmG convergence. 
\begin{thm} [\cite{Sturm06, AGS14b, AGMR15}] {\bf (Stability of the RCD$(K,\infty)$ condition)} \label{thm: Stab RCD}
\\
Let $\{\mathcal X_n\}_{n \in \N}$ be a sequence of RCD$(K,\infty)$ spaces.
If $\mathcal X_n$ converges to $\mathcal X_\infty$ in the pmG sense, then the limit space $\mathcal X_\infty$ is also an RCD$(K,\infty)$ space. 
\end{thm}

\subsubsection{$W_2$-gradient Flow of Relative Entropy} \label{subsub: RE}
In this subsection, following \cite{AGMR15, AGS14}, we recall the heat flow on the $L^2$-Wasserstein space  $(\mathcal P_2(X), W_2)$, which is constructed by the gradient flow of the relative entropy functional. We also recall the stability of the heat flows under the pmG convergence. 

 The descendent slope $|D^{-}{\rm Ent}_m|: \mathcal P_2(X) \to [-\infty, \infty]$ of the relative entropy ${\rm Ent}_m$ is defined as follows:
 \begin{align*}
 |D^-{\rm Ent}_m|(\mu)=
 \begin{cases} \dis
 \limsup_{W_2(\nu, \mu)\to 0}\frac{({\rm Ent}_m(\mu)-{\rm Ent}_n(\nu))^+}{W_2(\nu,\mu)}, &\text{if $\mu \in D({\rm Ent}_m)$},
  \vspace{3mm}\\\dis
 0, &\text{if $\mu$ is isolated in $\mathcal P_2(X)$},
 \vspace{3mm} \\ \dis
 +\infty, &\text{if $\mu \in \mathcal P_2(X)\setminus D({\rm Ent}_m)$}.
 \end{cases}
 \end{align*}
Here $(\cdot)^+$ denotes the positive part. 
 Let $\mathcal X=(X,d,m,\x)$ be a CD$(K,\infty)$ space and $\overline{\mu} \in D({\rm Ent}_m)$. A curve $\mu: [0,\infty) \to D({\rm Ent}_m) \subset \mathcal P_2(X)$ is said to be {\it the $W_2$-gradient flow of ${\rm Ent}_m$ starting at $\overline{\mu}$} if $\mu$ is locally absolutely continuous in $(\mathcal P_2(X), W_2)$ with $\mu_0=\overline{\mu}$ and 
 \begin{align*}
 {\rm Ent_m}(\mu_t)={\rm Ent_m}(\mu_s) + \frac{1}{2}\int_t^s|\dot{\mu}_r|^2dr+\frac{1}{2}\int_t^s|D^-{\rm Ent}_m|^2(\mu_r)dr, \quad 0 < \forall t < \forall s.
 \end{align*}
Under the CD$(K,\infty)$ condition, it is known that the gradient flow $\mu_t=\mathcal H_t\overline{\mu}$ of the relative entropy exists uniquely for any initial measure $\overline{\mu} \in \overline{D({\rm Ent}_m)}$ and for any $t \ge 0$ (\cite{AGMR15, AGS14}). Here $\overline{D({\rm Ent}_m)}$ means the closure of $D({\rm Ent}_m)$. We call $\{\mathcal H_t\}_{t \ge 0}$ {\it heat flow on $\mathcal P_2(X)$}.

\begin{thm}[Theorem 7.7 in \cite{GMS13}] \label{thm: CONHF} {\bf (Stability of heat flows)}
\\
Let $\{\mathcal X_n=(X_n,d_n,m_n,\x_n)\}_{n \in \N}$ be a sequence of  RCD$(K,\infty)$ spaces converging to $\mathcal X_\infty=(X_\infty, d_\infty, m_\infty, \x_\infty)$ in the pmG sense.  
If ${\mu}_n \in \mathcal P_2({\rm supp}[m_n]) \subset \mathcal P_2(X)$ converges to $\mu_\infty \in \mathcal P_2({\rm supp}[m_\infty])\subset \mathcal P_2(X)$ in the $W_2$-sense:
$$W_2({\iota_n}_\#{\mu}_n, {\iota_\infty}_\#{\mu}_\infty) \to 0, \quad n \to \infty,$$
then the solution $\mu_t^n=\mathcal H^n_t({\mu}_n)$ of the heat flow starting at ${\mu}_n$ converges to the limit one $\mu_t^\infty=\mathcal H^\infty_t({\mu}_\infty)$ in the $W_2$-sense:
$$W_2({\iota_n}_\#{\mu}^n_t, {\iota_\infty}_\#{\mu}^\infty_t) \to 0, \quad n \to \infty, \quad \forall t\ge0.$$
Here $\iota_n$ is an embedding $X_n \to X$ corresponding to the pmG convergence (see Definition \ref{prop: Dconv}).
\end{thm}

\subsubsection{$L^2$-gradient Flow of Cheeger Energy} \label{subsub: CE}
We now recall the $L^2$-gradient flow of Cheeger energy by Hilbertian theory of gradient flows (see e.g., \cite{AGS05}). We also recall the important fact that the heat flow in the previous section and the $L^2$-gradient flow of Cheeger energy in this section coincide under the CD$(K,\infty)$ condition.

For $f_0 \in L^2(X;m)$, there exists a locally Lipschitz map $t \mapsto f_t=H_tf_0 \in L^2(X;m)$ with $f_t \to f_0$ as $t \downarrow 0$ whose derivative satisfies 
\begin{align} \label{eq: GFC}
\frac{d}{dt}f_t \in -\partial^-{\sf Ch}(f_t), \quad \text{a.e.-}t > 0.
\end{align}
Here the subdifferential $\partial^-{\sf Ch}$ of convex analysis is the multi-valued operator in $L^2(X;m)$ defined at all elements of the domain of the Cheeger energy $f \in W^{1,2}(X,d,m)$ by the family of inequalities 
$$h \in \partial^-{\sf Ch}(f) \iff \int_Xh(g-f)dm \le {\sf Ch}(g)-{\sf Ch}(f), \quad \forall g \in L^2(X;m).$$
The map $H_t: f_0 \mapsto f_t$ is uniquely determined by \eqref{eq: GFC} and define a contraction semigroup (not necessarily linear) on $L^2(X;m)$. The flow $f_0 \mapsto f_t=H_tf$ is called {\it $L^2$-gradient flow of the Cheeger energy}, and the semigroup $\{H_t\}_{t \ge 0}$ is called {\it heat semigroup}. 

We recall that the $L^2$-gradient flows of Cheeger energies and the $W_2$-gradient flow of entropies are equivalent under the CD$(K,\infty)$ condition.
\begin{thm} \label{thm: CEE} \cite[Theorem 9.3]{AGS14}(see also \cite{AGMR15})
Let $\mathcal X=(X,d,m,\x)$ be a p.m.m.\ space satisfying the CD$(K,\infty)$ condition. If $\mu_0=f_0m \in \mathcal P_2(X)$ with $f_0 \in L^2(X;m)$, then 
$$\mathcal H_t(\mu_0)=(H_tf_0)m, \quad \forall t\ge 0.$$
\end{thm}

\subsubsection{$\mathrm{RCD^*}(K,N)$ Spaces} 
In this subsection, we recall the definition of the $\mathrm{RCD}^*(K,N)$ condition and several properties satisfied by $\mathrm{RCD}^*(K,N)$ spaces (see \cite{G15, AMS16, EKS15} for more details).

For each $\theta \in [0,\infty)$, we define the following functions
$$
\Theta_\kappa(\theta)=
\begin{cases} \dis
\frac{\sin(\sqrt{\kappa}\theta)}{\sqrt{\kappa}}, & \text{if} \quad \kappa>0,
\\ 
\dis
\theta,  &\text{if} \quad \kappa=0,
\\
\dis
\frac{\sinh(\sqrt{-\kappa}\theta)}{\sqrt{-\kappa}}, &\text{if} \quad \kappa<0,
\end{cases}
$$
We define the following functions: for $t \in [0,1]$, 
\begin{align*} \dis
\sigma_{\kappa}^{(t)}(\theta)=
\begin{cases} \dis \frac{\Theta_{\kappa}(t\theta)}{\Theta_{\kappa}(\theta)}, &\text{if} \quad \kappa \theta^2 \neq0 \quad \text{and} \quad \kappa \theta^2< \pi^2,
\\
t, &\text{if} \quad \kappa \theta^2=0,
\\
\dis
+\infty, &\text{if} \quad \kappa \theta^2 \ge \pi^2.
\end{cases}
\end{align*}
Let $\mathcal \mathbb P_\infty(X,d,m)$ be the subset of  $\mathcal P_2(X)$ consisting of $\mu$ which is absolutely continuous with respect to $m$ and has bounded support.
\begin{defn}[\cite{BS10, G15}] \normalfont ($\mathrm{CD}^*(K,N)$ and $\mathrm{RCD}^*(K,N)$)
\begin{description}
\item[(i)] A metric measure space $(X,d,m)$ is said to satisfy {\it the reduced curvature-dimension condition} $\mathrm{CD}^*(K,N)$ for $K, N \in \R$ with \replaced{$1<N<\infty$}{$N>1$} if,
for each pair $\mu_0=\rho_0m$ and  $\mu_1=\rho_1m$ in $\mathcal \mathbb P_\infty(X,d,m)$, there exist an optimal coupling $q$ of $\mu_0$ and $\mu_1$  and a geodesic $\mu_t=\rho_tm$ $(t \in [0,1])$ in $(\mathcal \mathbb P_\infty(X,d,m), W_2)$ connecting $\mu_0$ and $\mu_1$ so that, for all $t \in [0,1]$ and $N' \ge N$, we have
\begin{align*} 
\int\rho_t^{-\frac{1}{N'}} d\mu_t \ge  \int_{X \times X} & \Bigl[\sigma_{\frac{K}{N'}}^{(1-t)}(d(x_0,x_1))\rho_0^{-\frac{1}{N'}}(x_0)
+\sigma_{\frac{K}{N'}}^{(t)}(d(x_0,x_1))\rho_1^{-\frac{1}{N'}}(x_1)\Bigr]dq(x_0,x_1).\notag
\end{align*}
\item[(ii)] A metric measure space $(X,d,m)$ is said to satisfy {\it the Riemannian curvature-dimension condition} $\mathrm{RCD}^*(K,N)$ if the following two conditions hold:
\begin{description}
	\item[(ii-a)] $\mathrm{CD}^*(K,N)$
	\item[(ii-b)] the infinitesimal Hilbertianity, that is the Cheeger energy {\sf Ch} is a quadratic form: 
	\begin{align*}
	2{\sf Ch}(u)+2{\sf Ch}(v)={\sf Ch}(u+v) + {\sf Ch}(u-v), 
	\\
	\hfill \forall u,v \in W^{1,2}(X,d,m). \notag
	\end{align*}
\end{description}
\end{description}
\end{defn}
\begin{rem} \normalfont
The RCD$^*(K,N)$ condition is stronger than the RCD$(K,\infty)$ condition. \added{If $X$ is an RCD$^*(K,N)$ space, then $X$ is locally compact by the local volume doubling property according to Bishop--Gromov inequality \cite[Proposition 3.6]{EKS15} (see also \cite[{Corollary} 2.4]{Sturm06-1})}.
\end{rem}

The RCD$^*(K,N)$ condition is stable under the pmG convergence. 
\begin{thm} [\cite{G15}] {\bf (Stability of RCD$^*(K,N)$)}
\\
Let $\{\mathcal X_n\}_{n \in \N}$ be a sequence of RCD$^*(K,N)$ spaces. 
If $\mathcal X_n$ converges to $\mathcal X_\infty$ in the pmG sense, then $\mathcal X_\infty$ is also an  RCD$^*(K,N)$ space.  
\end{thm}

 \subsection{$L^2$-convergence of Heat Semigroups Under the PmG Convergence} \label{sec: MC}
In Gigli--Mondino--Savar\'e \cite{GMS13}, they introduced $L^2$-convergences on varying metric measure spaces and showed a convergence of heat semigroups in this sense under the pmG convergence of the underlying spaces with the RCD$(K,\infty)$ condition.  We recall their results briefly.

\begin{defn}$(${\rm See \cite[Definition 6.1]{GMS13}}$)$ \label{defn: Weak} \normalfont
Let $\{(X_n, d_n, m_n, \x_n)\}_{n \in \N}$ be a sequence of p.m.m. spaces. Assume that $(X_n,d_n,m_n, \x_n)$ converges to $(X_\infty, d_\infty, m_\infty, \x_\infty)$ in the pmG sense.  Let $(X,d)$ be a complete separable metric space and
$\iota_n: \mathrm{supp}[m_n] \to X$ be isometries as in Definition \ref{prop: Dconv}. We identify $(X_n, d_n, m_n)$ with $(\iota_n(X_n), d, {\iota_n}_\#m_n)$ and omit $\iota_n$. 
\begin{description}
	\item[(i)] We say that {\it $u_n \in L^2(X, m_n)$ converges weakly to $u_\infty \in L^2(X,m_\infty)$} if the following hold:
$$\sup_{n \in \N}\int|u_n	|^2\ dm_n < \infty \quad \text{and} \quad \int \phi u_n \ dm_n \to \int \phi u_\infty \ dm_\infty \quad \forall \phi \in C_{bs}(X),$$
whereby recall that $C_{bs}(X)$ denotes the set of bounded continuous functions with bounded support.
 \item[(ii)] We say that {\it $u_n \in L^2(X,m_n)$ converges strongly to $u_\infty \in L^2(X, m_\infty)$} if $u_n$ converges weakly to $u_\infty$ and the following holds:
$$\limsup_{n \to \infty}\int |u_n|^2 \ dm_n \le \int|u_\infty|^2\ dm_\infty.$$
\end{description}
\end{defn}
Let $\{H^n_t\}_{t\ge 0}$ be the $L^2(X,m_n)$-semigroup corresponding to the Cheeger energy ${\sf Ch}_n$. Then the following theorem states that $\{H^n_t\}_{t\ge 0}$ convergence strongly in $L^2$ under the pmG convergence of the underlying spaces. 
\begin{thm} $(${\rm See \cite[Theorem 6.11]{GMS13}}$)$ \label{thm: Mosco of SG}
Let $\{(X_n, d_n, m_n, \x_n)\}_{n \in \N}$ be a sequence of p.m.m.\ spaces satisfying the $\mathrm{RCD}(K,\infty)$ for all $n \in \N$. Then, for any $u_n \in L^2(X,m_n)$ converging strongly to $u_\infty \in L^2(X,m_\infty)$, we have, for any $t>0$
$$H_t^nu_n \text{ converges strongly to } H_t^\infty u_\infty \text{ in the sense of Definition \ref{defn: Weak}}.$$
\end{thm}
Note that, in \cite[Theorem 6.11]{GMS13}, the above Theorem \ref{thm: Mosco of SG} was stated without the condition of the infinitesimal Hilbertian. 

\section{Brownian Motion on RCD spaces} \label{sec: BM}
\subsection{Brownian Motions on RCD$(K,\infty)$ Spaces}\label{subse: BM1}
Let $(X,d,m)$ satisfy the $\mathrm{RCD}(K,\infty)$ condition.
Let $\delta_x$ denote the unit mass at $x \in X$, and define a kernel $p(t,x,dy)$ by the action of the heat flow (see Subsection \ref{subsub: RE})
\begin{align*}
p(t,x,dy):=\mathcal H_t(\delta_x) \quad \forall t>0, x \in X.
\end{align*}
Then we have that (see \cite{AGS14b, AGMR15})
$$\text{$p(t,x,dy)$ is absolutely continuous with respect to $m$ for any $t>0$,}$$ 
and we denote the density by $p(t,x,y)$. By \cite[Theorem 6.1]{AGS14b} and \cite{AGMR15} (for the case of $\sigma$-finite reference measures), the density $p(t,x,y)$ is symmetric \replaced{in $x$ and $y$}{(i.e. $p(t,x,y)=p(t,y,x)$ for {any $x, y \in {\rm supp}[m]$})}, and satisfies the Chapman--Kolmogorov formula. Moreover, the following action of semigroup $\{P_t\}_{t \ge 0}$
\begin{align} \label{eq: HSBM}
P_tf(x):=\int_X f(y) d \mathcal H_t(\delta_x)(dy)
\end{align}
is a version of the linear heat semigroup $\{H_t\}_{t \ge 0}$ defined as the gradient flow of the Cheeger energy ${\sf Ch}$ (see Subsection \ref{subsub: CE}) for any $f \in L^2(X;m)$. Furthermore $P_t$ is an extension of $H_t$ to a continuous contraction semigroup in $L^1(X;m)$ which is {\it point-wise everywhere defined on ${\rm supp}[m]$} if $f \in L^\infty(X;m)$ since $P_tf$ becomes Lipschitz continuous on ${\rm supp}[m]$ whenever $f \in L^\infty(X;m)$ (see \cite[Theorem 6.5]{AGS14b} and \cite[Theorem 7.3]{AGMR15}). 
We call $p(t,x,dy)$ and $p(t,x,y)$ {\it the heat kernel } and {\it the heat kernel density}, respectively. 
By the Kolmogorov extension theorem, we can construct a family of probability measures $\{\mathbb P^x\}_{x \in X} $ on $X^{[0,\infty)}$ and a system of Markov processes $(\{\mathbb P^x\}_{x \in X}, \{B_t\}_{t \ge 0})$ on $X$ with respect to $p(t,x,dy)$. 

On the other hand, we can define a Dirichlet form (i.e., a symmetric closed Markovian bilinear form) $(\mathcal E, \mathcal F)$ induced by the Cheeger energy ${\sf Ch}$ as follows:
\begin{align*} 
\mathcal E(u,v)=\frac{1}{4}({\sf Ch}(u+v)-{\sf Ch}(u-v)), \quad u,v \in \mathcal F=W^{1,2}(X,d,m).
\end{align*}
By {\cite[Lemma 6.7]{AGS14b}} (see \cite[Theorem 7.2]{AGMR15} for $\sigma$-finite reference measures), the form $(\mathcal E, \mathcal F)$ becomes a quasi-regular conservative strongly-local symmetric Dirichlet form below. 
See \cite[Proposition 4.11]{AGS14b} for the strong locality, and the conservativeness follows from the volume growth estimate \eqref{asmp: basic1} and Sturm's conservativeness test \cite[Theorem 4]{Sturm94}. 

Therefore, by \cite[Theorem IV 3.5, V1.5]{MR92}, there exists a family of probability measures $\{\mathbb Q^x\}_{x \in X}$ on $C([0,\infty);X)$ and a system $(\{\mathbb Q^x\}_{x \in X}, \{B'_t\}_{t \ge 0})$ of conservative diffusion processes so that 
$$\mathbb E_{\Q}^x (f(B'_t))=P_tf, \quad \forall f \in L^2(X;m)\cap \mathcal B_b(X), \quad \forall t \ge 0, \quad \forall x \in X \setminus \mathcal N.$$
Here $\mathbb E_{\Q}^x$ denotes the expectation with respect to $\mathbb Q^x$ and $\mathcal N$ is a set of zero-capacity with respect to $(\E, \F)$. By {\it conservative diffusion process}, we mean that $(\mathbb Q^x, \{B'_t\}_{t \ge 0})$ is a strong Markov process whose sample path is continuous, which means $B'_{\cdot} \in C([0,\infty); X)$ $\mathbb Q^x$-almost surely for every $x \in X$. The systems of Markov processes corresponding to $(\E, \F)$ are unique up to \replaced{zero}{null}-capacity sets with respect to starting points $x$. Namely, if there is another system of diffusion processes $(\{\mathbb R^x\}_{x \in X}, \{S_t\}_{t \ge 0})$ corresponding to $(\E, \F)$, then the laws of $(\mathbb Q^x, \{B'_t\}_{t \ge 0})$ and $(\mathbb R^x, \{S_t\}_{t \ge 0})$ coincide for every $x \in X \setminus \mathcal N$ with some set $\mathcal N$ of zero-capacity. Note that, if $\{P_t\}_{t \ge 0}$ is a Feller semigroup, then $\mathcal N$ can be taken as an empty set $\emptyset$, and the diffusion process can be defined uniquely with respect to every starting point $x$. See, e.g., \cite[Chapter 7, A.2]{FOT11} and \cite[Chapter IV]{MR92} for more comprehensive accounts.
\deleted{By considering that these two processes are equivalent, we call it {\it a system of Brownian motions}.} \deleted{(note that $(\{\mathbb P^x\}_{x \in X}, \{X_t\}_{t \ge 0})$ is defined at every starting point by the Kolmogorov extension theorem).}

Let ${\mathbb P}^x_*$ denote the outer measure of $\mathbb P^x$ on all subsets of $X^{[0,\infty)}$. By the same argument of \cite[\added{Proof of }(c) in Theorem 1.2]{ASZ09}, we have that 
\begin{align} \label{conserv}
\mathbb P^x_*\bigl(C((0, \infty);X)\bigr)=1\quad  \text{for {\it every} $x \in X$ (not only quasi-every $x \in X$).}
\end{align}
\deleted{Note that by the conservativeness \added{and the strong locality} of the Dirichlet form $(\mathcal E, \mathcal F)$, we know that $\mathbb P^x_*(C([0, \infty)))=1$ for {\it quasi-every $x \in X$}. However, the property \eqref{conserv} is not necessarily true for general conservative quasi-regular strongly local Dirichlet form, and}\replaced{T}{t}his property is due to the absolute continuity of the heat kernel $p(t,x,dy)$ with respect to $m$ for any $t>0$.
Note that, by the conservativeness \added{and the strong locality} of the Dirichlet form $(\mathcal E, \mathcal F)$, we know that $\mathbb Q^x(C([0, \infty);X))=1$ for {\it quasi-every $x \in X$} (this holds also for $\mathbb P_*^x$). We write simply $\mathbb P^x$ for ${\mathbb P}^x_*$. The systems of Markov processes $(\{\mathbb P^x\}_{x \in X}, \{B_t\}_{t \ge 0})$ and $(\{\mathbb Q^x\}_{x \in X}, \{B'_t\}_{t \ge 0})$ coincide except on zero-capacity sets. \replaced{In this paper}{To remove the ambiguity of starting points}, we adopt $(\{\mathbb P^x\}_{x \in X}, \{B_t\}_{t \ge 0})$ for representing a system of Brownian motions. 
\begin{rem} \normalfont \label{rem: BMN}
The diffusion process defined above is conventionally called {\it Brownian motion} (\cite{AGS14b}), but this may indicate other diffusion processes than the standard Brownian motion in some situations. For instance, when we take $(X,d,m)=(\R^d, \|\cdot\|_2, \frac{1}{(2 \pi)^{d/2}}\exp\{-\frac{1}{2}\|x\|_2^2\}dx)$ whereby $dx$ denotes the Lebesgue measure, and $\|\cdot\|_2$ denotes the Euclidean distance. Then $(X,d,m)$ satisfies RCD$(0,\infty)$ and the diffusion induced by the Cheeger energy coincides with what is known as the Ornstein-Uhlenbeck process, which is different from the standard Brownian motion on $\R^d$. 
\end{rem}

\subsection{Brownian Motions on RCD$^*(K,N)$ Spaces} \label{subse: BM}
In this subsection, \replaced{we show the Feller property of the heat semigroup on ${\rm RCD}^*(K,N)$ spaces.
}{we state several properties of the Brownian motions on ${\rm RCD}^*(K,N)$ spaces with bounded diameters, which do not necessarily hold only for the RCD$(K,\infty)$ condition.
The main points are that the Brownian motion under Assumption \ref{asmp: 1} satisfies a uniform Gaussian estimate and becomes a Feller process, which implies the uniqueness of Brownian motions with respect to starting points. }

\begin{prop}  \label{prop: Feller_BM}
Under \replaced{the RCD$^*(K,N)$ condition}{Assumption \ref{asmp: 1}}, 
\added{the heat semigroup $\{H_t\}_{t \ge 0}$ has a Feller modification. That is, there exists a}\deleted{modification} \added{semigroup $\{P_t\}_{t \ge 0}$ so that $P_tf=H_tf$ $m$-a.e.\ for any $f \in L^2(X,m)$ and any $t>0$ and the following conditions hold:} \deleted{there exists a modification semigroup $\{P_t\}_{t \ge 0}$ of the heat semigroup $\{H_t\}_{t \ge 0}$ so that $\{P_t\}_{t \ge 0}$ becomes a Feller semigroup, i.e.,} 
\begin{description}
\item[(F-1)] For any $f \in C_{{\infty}}(X)$, $P_tf \in C_{{\infty}}(X)$ for any $t>0$.
\item[(F-2)] For any $f \in C_{{\infty}}(X)$, $\|P_tf-f\|_\infty \to 0 \quad t \downarrow 0$.
\end{description}
\end{prop}
\begin{rem} \normalfont
The following proof is the result of a private communication with Prof.\ Kazuhiro Kuwae.
Although the proof might be already known in some literature, we could not find good references and we give the proof for the sake of reader's convenience. 
\end{rem}
\proof \added{By} {\cite[(iii) in Theorem 6.1]{AGS14b},} \added{there exists a semigroup $\{P_t\}_{t \ge 0}$ which is a modification of $\{H_t\}_{t \ge 0}$ so that $P_tf \in {\rm Lip}_{b}(X)$ if $f \in L^{\infty}(X,m)$. Before checking (F-1) and (F-2), we first give a heat kernel estimate.}
{By \cite[Theorem 1.2]{JLZ15}, we have the following Gaussian heat kernel estimate:
 there exist positive constants $C_i=C_i(N,K)$ for $i=1,2,3$ depending only on $N,K$ so that}
 \begin{align} \label{GE1}
  p(t,x,y) \le \frac{C_1}{m(B_{\sqrt{t}}(y))}\exp\Bigl\{-C_2\frac{d(x,y)^2}{t}{-C_3t}\Bigr\},
 \end{align}
for all $x, y \in X$ and\deleted{$\text{$0 < t <D^2$}$} \added{$0<t$}. Here the heat kernel $p(t,x,y)$ means the integral kernel of the heat semigroup $P_tf(x)=\int_X fp(t,x,y)m(dy)$ for $t>0$. 

\added{We now show condition (F-1). We already know $P_tf \in C_b(X)$, so it suffices to show that $P_tf$ vanishes at infinity for $f \in C_\infty(X)$ and $t>0$, which is an easy consequence of \eqref{GE1} as follows: we may assume that $f$ is compactly supported since every element in $C_\infty(X)$ can be approximated by elements in $C_0(X)$ with respect to the uniform norm. Let $K \supset {\rm supp}[f]$ be a compact set. By \eqref{GE1} and $\inf_{y \in K}m(B_{\sqrt{t}}(y))>0$ (by the lower semi-continuity of $m(B_r(x))$ in $x$), we see that, for any $\e>0$, there exists a compact set $K' \subset X$ so that}
{\begin{align*}
|P_tf(x)| &\le \int_{K} p(t,x,y)|f(y)|m(dy)
\\
&< \|f\|_\infty \int_{K} \frac{C_1}{m(B_{\sqrt{t}}(y))}\exp\Bigl\{-C_2\frac{d(x,y)^2}{t}{-C_3t}\Bigr\}m(dy)
\\
&< \frac{C_1\|f\|_\infty}{\inf_{y \in K}m(B_{\sqrt{t}}(y))} \int_{K} \exp\Bigl\{-C_2\frac{d(x,y)^2}{t}{-C_3t}\Bigr\}m(dy)
\\
& \le \e \quad (\forall x \in X \setminus K'). 
\end{align*}
Thus we \added{have} proved  (F-1).

Now we prove (F-2). We may assume $f \in C_0(X)$. 
Let $K \supset {\rm supp}[f]$ be a compact set. For given $\e>0$, take $\delta>0$ so that $|{f}(x)-{f}(y)|<\e$ whenever $d(x,y)<\delta$ in $x,y \in K$. By the Gaussian estimate \eqref{GE1}, we can choose a positive number $T$ so that $p(t,x,y)<\e$ for any $0<t<T$, and for any $x\in X$ and $y \in K$ satisfying $d(x,y) \ge \delta$.
Then we have that, for any $x \in X$
\begin{align*}
&|P_t{f}(x)-{f}(x)|=\Bigl|\int_{K}p(t,x,y){f}(y)m(dy)-{f}(x)\Bigr| 
\\
&\le \int_{K}p(t,x,y)\Bigl|{f}(y)-{f}(x)\Bigr|m(dy)
\\
&=\int_{B_\delta(x) \cap K}p(t,x,y)|{f}(y)-{f}(x)|m(dy) + \int_{(B_\delta(x))^c\cap K)}p(t,x,y)|{f}(y)-{f}(x)|m(dy)
\\
&\le\  \e + 2\e\|f\|_\infty.
\end{align*}
Thus we have shown that (F-2) holds. 
\qed

 \section{Examples} \label{sec: exa}
In this section, several specific examples satisfying Assumption \ref{asmp: 2} or \replaced{the assumption in Theorem \ref{thm: mthm1}}{Assumption \ref{asmp: 1}} are given. 
In the first subsection, we explain weighted Riemannian manifolds whose weighted Ricci curvature is bounded below, and \replaced{their}{its} pmG limit spaces. 
In the second subsection, we explain Alexandrov spaces, which \replaced{are}{is} \replaced{a generalization of}{roughly speaking,} the \added{lower} sectional curvature\deleted{is} bound\deleted{ed flow below} \added{to metric spaces}. In the third subsection, we give Hilbert spaces with log-concave probability measures. 
\subsection{Weighted Riemanniam Manifolds and pmG Limit Spaces} 
Let $\{(M_n, g_n, w_n, \x_n)\}_{n \in \N}$ be a sequence of pointed complete and connected weighted $N$-dimensional Riemannian manifolds whose weight satisfies $w_n=e^{-V_n}$ for a twice continuously differentiable function $V_n \in C^2(M_n)$. We write the corresponding pointed metric measure space $\mathcal M_n=(M_n, d_{g_{n}}, w_n{\rm Vol}_n, \x_n)$ whereby $d_{g_{n}}$ denotes the distance function associated with the Riemannian metric $g_n$; ${\rm Vol}_n$ denotes the Riemannian volume measure; and $\x_n \in M_n$ is a fixed point. 
 Let the weighted Ricci curvature ${\rm Ric}_{\mathcal M_n}$ of $\mathcal M_n$ be bounded from below by $K$: there exists $K \in \R$ so that 
$${\rm Ric}_{\mathcal M_n}={\rm Ric}_{g}+\nabla^2V_n \ge Kg_n,$$
whereby ${\rm Ric}_{g_n}$ means the Ricci curvature of $(M_n,g_n)$ and $\nabla^2$ means the Hessian. Then $\mathcal M_n$ satisfies RCD$(K,\infty)$ spaces (\cite{RS05, Sturm06}). Even when $V_n: M_n \to \R$ is not in $C^2(M_n)$, if ${\rm Ric}_{g_n} \ge K$ and 
$$\text{$V_n: M_n \to \R$ is $K'$-convex (see \cite{Sturm06})},$$
 then $\mathcal M_n$ satisfies RCD$(K+K',\infty)$. 
 If, moreover, 
 $$\text{$V_n: M_n \to \R$ is $(K',N')$-convex (see \cite{EKS15}),}$$
  then $\mathcal M_n$ satisfies RCD$^*(K+K',N+N')$.
The Brownian motion on $M_n$ is a Markov process whose infinitesimal generator $A_n$ is 
$$A_n=\frac{1}{2}\Delta_{M_n}-\langle \nabla V_n, \nabla \rangle,$$
whereby $\Delta_{M_n}$ is the Laplace-Beltrami operator on $M_n$.

If $\mathcal M_n$ satisfying RCD$(K,\infty)$ (or, RCD$^*(K,N)$) converges to a metric measure space $\mathcal M_\infty$ in pmG sense, then the limit space $\mathcal M_\infty$ satisfies RCD$(K,\infty)$ (or, RCD$^*(K,N)$), respectively (see \cite{GMS13, EKS15}). Thus we can apply our main results\deleted{(Theorem \ref{thm: mthm2}, or Theorem \ref{thm: mthm2-1} for the case of RCD$^*(K,N)$) to Brownian motions on these spaces} and\deleted{we can} obtain the weak convergence of the Brownian motions.  

We have various singular examples appearing as the limit space. See e.g., \cite[Example 8]{CC97}. We give one of the simplest examples included in this framework. 
 \begin{exa} \normalfont {\bf (Collapsing: Torus $\to$ Circle)} \vspace{-1mm} \label{ex: col}
 \\
 Let $\mathbb S^1 \subset \R^2$ be the unit circle. Let $d_{\mathbb S^1}$ be the \replaced{shortest path}{intrinsic} distance on $\mathbb S^1$, that is, the distance between $x$ and $y$ is defined by the infimum over lengths of geodesics on $\mathbb S^1$ connecting $x$ and $y$. 
Let $$\overline{H}_{\mathbb S^1}:=\frac{1}{H_{\mathbb S^1}(S^1)}H_{\mathbb S^1}$$ be the normalized Hausdorff measure on $(\mathbb S^1, d_{\mathbb S^1})$.
	Let $\mathbb T_n=\mathbb S^1 \times \mathbb S^1$ be a two-dimensional flat torus with a metric $d_n=d_{\mathbb S^1}\otimes \frac{1}{n}d_{\mathbb S^1}$ and the normalized Hausdorff measure $\overline{H}_n$ on $(\mathbb T_n, d_n)$, whereby
	$$d_{\mathbb S^1}\otimes \frac{1}{n}d_{\mathbb S^1}((x_1,y_1),(x_2,y_2)):=\sqrt{d^2_{\mathbb S^1}(x_1,x_2)+\frac{1}{n^2}d^2_{\mathbb S^1}(y_1,y_2)}.$$
	Then $(\mathbb T_n, d_n, \overline{H}_n)$ satisfies the $\mathrm{RCD}^*(0, 2)$ for any $n \in \N$ and converges to $(\mathbb S^1, d_{\mathbb S^1}, \overline{H}_{\mathbb S^1})$ in the measured Gromov sense. Thus we can apply our result (Theorem \ref{thm: mthm1}) and the weak convergence of the Brownian motions\deleted{on $M_n$} is equivalent to the pmG convergence of \replaced{the underlying spaces}{$M_n$}.
	\begin{figure}[htbp]
\begin{center}\includegraphics[width=10cm, bb=0 0 842 295]{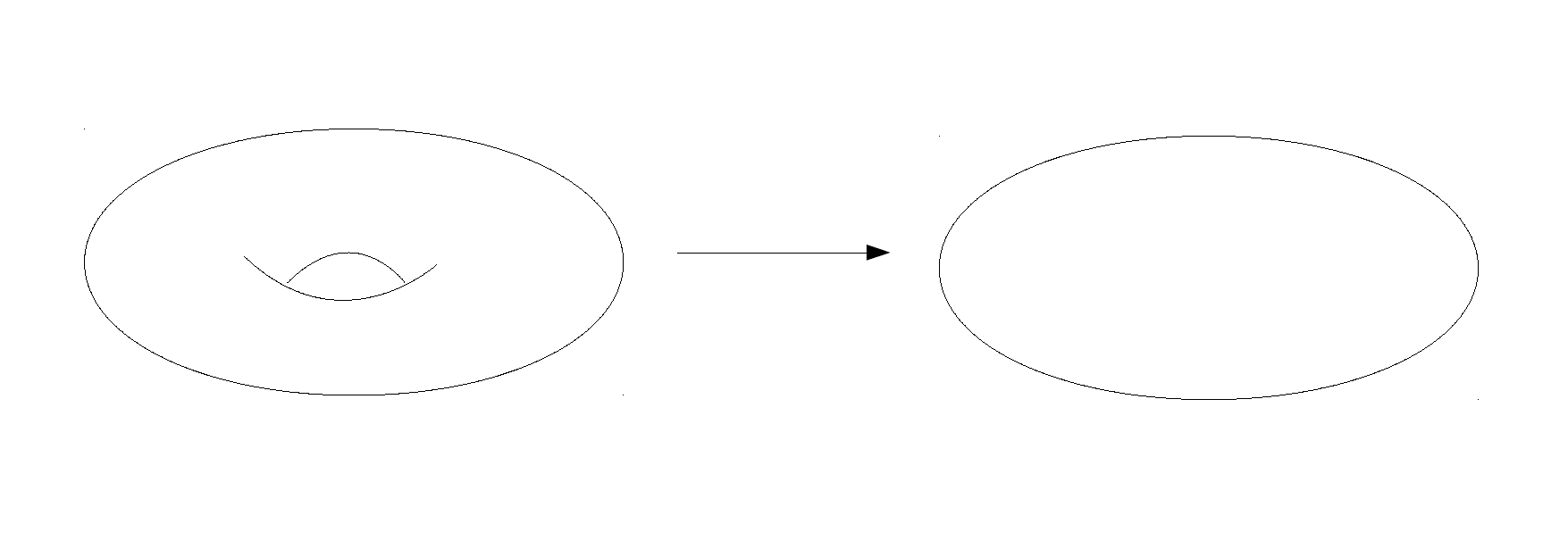}
\caption{Tori Converge to a Circle.}
\label{picture.1}
\end{center}
\end{figure}
\end{exa}

\subsection{Alexandrov Spaces}
We explain Alexandrov spaces, which are \added{a} generalization of lower bounds of sectional curvatures \replaced{to}{for} metric spaces.  We refer \replaced{the}{a} reader to \cite{BBI01} for basic theory of Alexandrov spaces. 
Let $(X,d)$ be a locally compact length space. For a triple of points $p,q,r \in X$, a geodesic triangle $\triangle pqr$ is a triplet of geodesics joining each two points. Let $M^N(K)$ be the $N$-dimensional complete simply connected space of constant sectional curvature $K$. For a geodesic triangle $\triangle pqr$, we denote by $\triangle \tilde{p} \tilde{q} \tilde{r}$ a geodesic triangle in $M^2(K)$ whose corresponding edges have the same lengths as $\triangle pqr$.

A locally compact length space $(X,d)$ is said to be {\it an Alexandrov space with ${\sf Curv} \ge K$} if for every point $x \in X$, there exists an open set $U_x$ including $x$ so that for every geodesic triangle $\triangle pqr$ whose edges are totally included in $U_x$, the corresponding geodesic triangle $\triangle \tilde{p} \tilde{q} \tilde{r}$ satisfies the following condition: for every point $z \in qr$ and $\tilde{z} \in \tilde{q}\tilde{r}$ with $d(q,z)=d(\tilde{q},\tilde{z})$, we have 
$$d(p,z) \ge d(\tilde{p}, \tilde{z}).$$
If we consider a complete $N$-dimensional Riemannian manifold $(M,g)$, then $(M,g)$ is an Alexandrov space with ${\sf Curv} \ge K$ if and only if ${\rm sec}(M) \ge K$, whereby ${\rm sec}(M)$ means the sectional curvature of $M$.

 Let $\mathcal X=(X,d, H)$ be an $N$-dimensional Alexandrov space with ${\sf Curv} \ge K$ and $H$ be the Hausdorff measure (see e.g., \cite{BBI01} for details). 
According to \cite{Pet11, ZZ10}, $\mathcal X$ satisfies $\mathrm{CD}^*((N-1)K,N)$. Moreover, as was shown in \cite{KMS01}, $\mathcal X$ satisfies the infinitesimal Hilbertian condition, and as a result, $\mathcal X$ satisfies $\mathrm{RCD}^*((N-1)K,N)$. Thus we can apply our results (Theorem \ref{thm: mthm2}, \ref{thm: mthm2-0}) and if a sequence of pointed Alexandrov spaces $X_n$ with ${\sf Curv} \ge K$ converges to the limit space $X_\infty$ in the pmG sense, then the Brownian motions on $X_n$ converge weakly to the limit Brownian motion on $X_\infty$.  
We give several examples.
\begin{exa} \normalfont
 {\bf (Cone $\to$ Interval)} \ \\
 Let $X_n \subset \R^3$ be a cone defined by $X_n=\{(x,y,z) \in \R^3: y^2+z^2=\frac{1}{n}x, \ 0 \le x <1\} \cup \{(x,y,z)\in \R^3: y^2+z^2=\frac{1}{n},\ x=1\}$. 
	Let $d_n$ be the \replaced{shortest path}{restriction of the Euclidean} distance on $X_n$ and $\overline{H}_n$ be the normalized Hausdorff measure on $X_n$.
	Then $(X_n,d_n, \overline{H}_n)$ satisfies $\mathrm{RCD}^*(0,2)$ and converges to $([0,1], |\cdot|, m)$ in the measured Gromov sense, whereby $m$ is \replaced{a}{some} measure on $[0,1]$. 
	 Thus we can apply our result (Theorem  \ref{thm: mthm1}) and the weak convergence of the Brownian motions \deleted{on $M_n$} is equivalent to the pmG convergence of \replaced{the underlying spaces}{$M_n$}.
	\begin{figure}[htbp]
\begin{center}\includegraphics[width=10cm, bb=10 10 842 357]{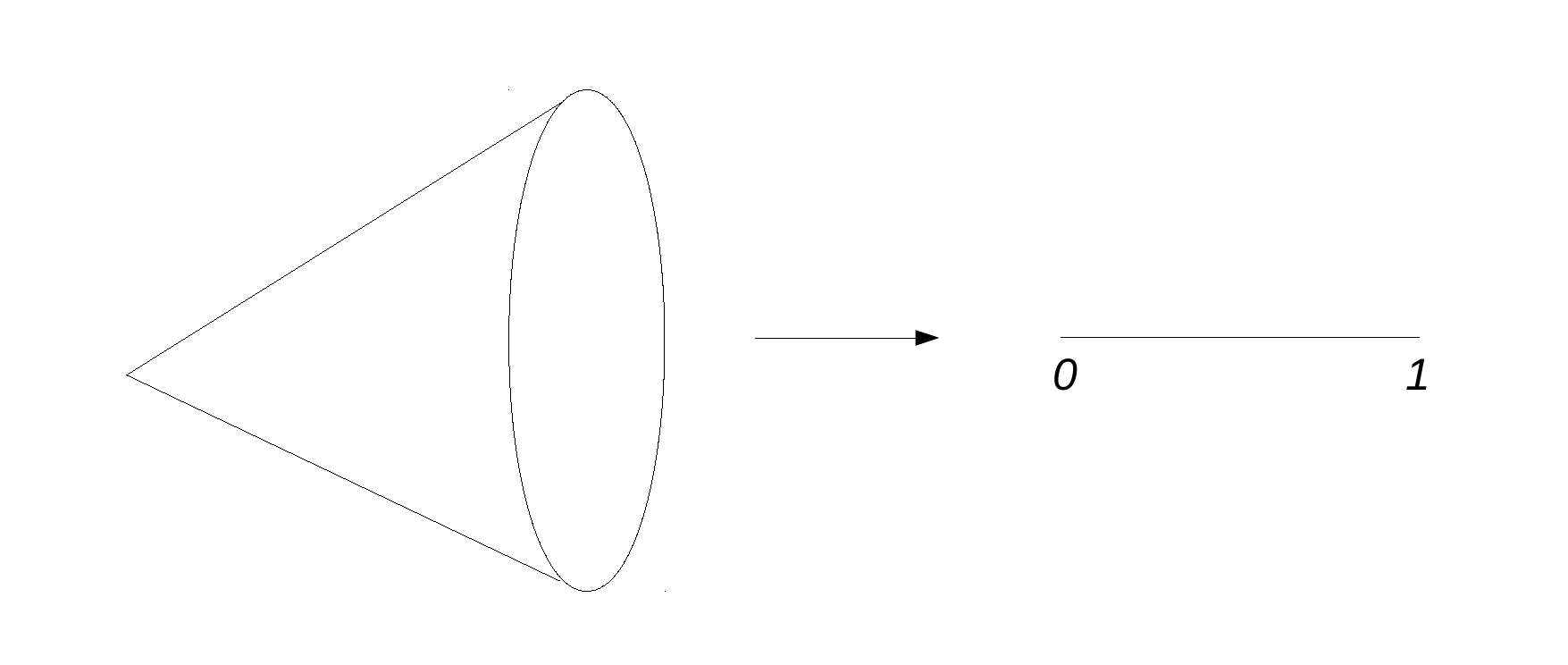}
\caption{Cones Converge to an Interval.}
\label{picture.2}
\end{center}
\vspace{-5mm}
\end{figure}
 \end{exa}

As a second example, we give a sequence of polygons made by the barycentric subdivision. The limit space has dense singularities. 
 \begin{exa}\normalfont {\bf (Dense Singularities \cite[p.\ 632, Examples.\ (2)]{OS94})}\  \label{exa: dens} \\
 Let $X=(X, d)$ be a polyhedron in $\R^3$ with the \replaced{shortest path}{geodesic} metric $d$ \added{on $X$}. Then we can check whether $X$ is an Alexandrov space with {\sf Curv}$\ge 0$, which is also an RCD$^*(0,2)$ space.
For any vertex $p \in X$, let $\angle(X, p)$ denote the sum of all inner angles at $p$ of faces $T$'s such that $p$ is a vertex of $T$. 

Now we construct a sequence of polyhedra $\{M_n\}_{n \in \N}$ inductively. Let $M_1$ be a tetrahedron in $\R^3$ with the barycenter $o$. 
Let $M_n$ be defined. Then we define $M_{n+1}$ as follows: Take a monotone decreasing sequence $\{\e_n\}_{n \in \N}$ \replaced{so}{such} that $\e_n \to 0$ as $n\to \infty$ with $0 < \e_n<1$ and $\e:=\Pi_{n=1}^\infty(1-\e_n)>0$. We take the barycentric subdivision of $M_n$. Keep the original vertices in $M_n$ \added{in the same positions} and move the new vertices generated by the barycentric subdivision outward along rays emanating from $o$ so small that, for the new polyhedra $M_{n+1}$ generated by the new and original vertices, we have 
 $$2\pi-\angle(M_{n+1},p) \ge (1-\e_n)(2\pi-\angle(M_n,p)),$$
 for any vertex $p \in M_n$. See \cite[p.\ 632, Examples.\ (2)]{OS94} for more details.
 
 \added{Let $d_n$ and $H_n$ be the shortest path distance and the Hausdorf measure on $M_n$.} Then there exists the Hausdorff-limit of $M_n=(M_n,d_n)$, denoted by $M_\infty$. The limit space $M_{\infty}$ is a two-dimensional Alexandorv space with nonnegative curvature.
 In particular, $(M_n,d_n,H_n)$ \deleted{satisfies the $\mathrm{RCD}^*(0, 2)$ and}converges to $(M_\infty, d_\infty, H_\infty)$ in the measured Gromov sense.  
The limit space $M_\infty$ also satisfies the $\mathrm{RCD}^*(0, 2)$ by the stability of RCD$^*(K,N)$ spaces under the measured Gromov convergence (see \cite{EKS15}). The set of singular points \added{in $M_\infty$} is dense\deleted{in $M_\infty$} (see \cite{OS94}). \replaced{Since each diameter of $M_n$ is obviously uniformly bounded by the construction,}{Thus} we can apply our result (Theorem  \ref{thm: mthm1}) and the weak convergence of the Brownian motions\deleted{on $M_n$} is equivalent to the pmG convergence of \replaced{the underlying spaces}{$M_n$}.
\begin{figure}[htbp] \vspace{0mm}
\begin{center}\includegraphics[width=10cm, bb=10 10 842 388]{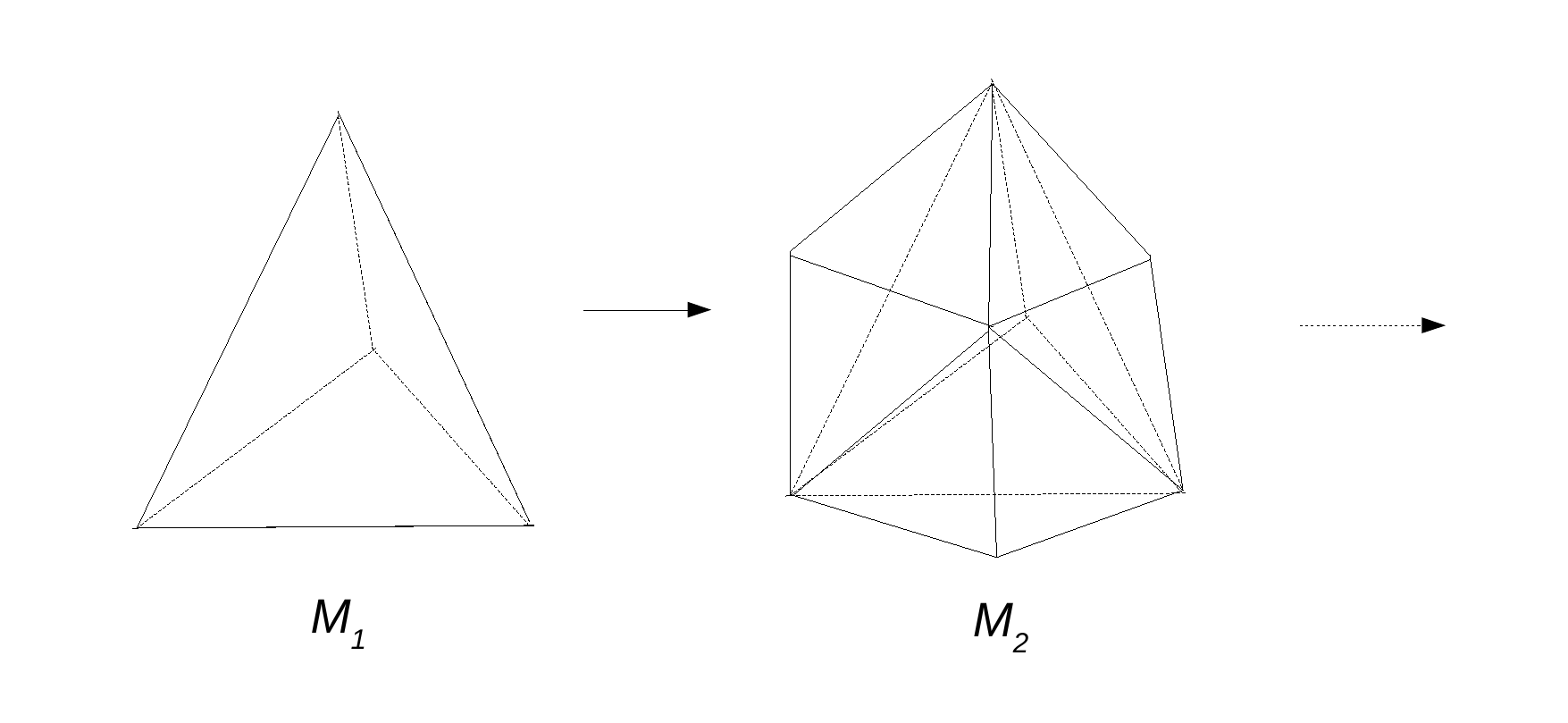}
\vspace{-6mm}
\caption{Polyhedra Generated by Barycentric Subdivision.}
\label{picture.3}
\end{center}
\vspace{0mm}
\end{figure}

 \end{exa}

 \subsection{Hilbert Space with Log-concave Measures}
 In this subsection, we give a specific class of RCD$(0,\infty)$ spaces, which is a Hilbert space with log-concave measures. This subsection follows \cite{ASZ09}.
 
 Let $H$ be a separable Hilbert space, which would be a finite- or infinite-dimensional space, with an inner product $\langle \cdot, \cdot \rangle$ and \replaced{the corresponding}{a} norm $\|\cdot\|$. A Borel probability measure $\gamma$ on $H$ satisfies {\it log-concave condition} if, for all pairs of open subsets $A, B \subset H$, it holds that 
$$\log{\gamma}((1-t)A+tB) \ge (1-t)\log \gamma(A)+t\log \gamma(B), \quad \forall t \in [0,1].$$
Let $K={\rm supp}[\gamma]$ and $A=A(\gamma)$ be the smallest closed linear subspace containing $K$. We write canonically 
$$A=H_0+h_0, \quad h_0 \in K, \quad \|h_0\| \le \|k\|, \quad \forall k \in K,$$
so that $h_0$ is the element of the minimal norm in $K$ and $H_0$ is a closed linear subspace in $H$.

Let $C_b^1(A)$ be the set of all $\Phi: A \to \R$ which are bounded, continuous and Fr\'echet differentiable with a bounded continuous gradient $\nabla \Phi: A \to H_0$. Then, according to \cite[Theorem 1.2]{ASZ09}, the following bilinear form becomes closable and the closed form becomes a symmetric quasi-regular Dirichlet form $\mathcal E=\mathcal E_{\|\cdot\|, \gamma}$:
\begin{align} \label{defn: SDE}
\mathcal E(u,v)=\int_K \langle \nabla u, \nabla v \rangle_{H_0}d\gamma, \quad u,v \in \mathcal F:=\overline{C_b^1(A)}^{\sqrt{\mathcal E+\|\cdot\|^2_2}}.
\end{align}
In \cite{ASZ09}, the corresponding semigroup $\{P_t\}_{t \ge 0}$ associated with $(\mathcal E, \mathcal F)$ satisfies EVI$_0$ property, which is equivalent to the RCD$(0,\infty)$ condition of $(H, \|\cdot\|, \gamma)$ according to \cite{AGS14b}. Let $\{H_n=(H,\|\cdot\|_n, \gamma_n, \x_n)\}_{n \in \N}$ be a sequence of pointed Hilbert spaces with log-concave probability measures satisfying the above conditions. Then 
the weak convergence of the Brownian motions on $H_n$ to that on $H_\infty$ follows from the pmG convergence of the underlying spaces $H_n$ to $H_\infty$ (Theorem \ref{thm: mthm2}, \ref{thm: mthm2-0} and \cite[Theorem 1.5]{ASZ09}).

Various infinite dimensional examples are included in the framework of Hilbert spaces with log-concave probability measures. For instance, all measures $\gamma$ of the following form satisfies the log-concave condition: let $dx$ be the Lebsgue measure on $\R^N$ and 
\begin{align*} 
\gamma=\frac{1}{Z}e^{-V}dx, \quad \text{whereby}\  V: H=\R^N \to \R\ \text{convex and} \ Z=\int_{\R^N}e^{-V}dx<+\infty,
\end{align*}
such as all Gaussian measures and all Gibbs measures on on a finite lattice with a convex Hamiltonian. See \cite[Section 1.2]{ASZ09} for various infinite-dimensional literatures related to stochastic partial differential equations. We give several finite-dimensional examples.
  \begin{exa}[\cite{ASZ09}] \label{ex: GUS}\normalfont
   We \replaced{explain}{give} several examples associated with stochastic differential equations (SDE). \replaced{The first one is}{We first give} SDEs on the Euclidean space $\R^N$ with variable potentials. \replaced{The second one is}{We secondly give} SDEs on variable convex domains in $\R^N$ with variable potentials. 
 \begin{enumerate}
 \item[(a)]{\bf (SDE with \replaced{V}{v}ariable \replaced{C}{c}onvex \replaced{P}{p}otentials)}
 Let $H=\R^N$ with $1<N<\infty$. Let $\{V_n: \R^N \to \R\}_{n \in \N}$ be a sequence of convex functionals with a Lipschitz continuous gradient $\nabla V_n: \R^N \to \R^N$ and $\int_{\R^N}e^{-V_n}dx<\infty$.
 Take
 \begin{align*}
\gamma_n=\frac{1}{Z_n}e^{-V_n}dx, \quad \text{whereby}\  Z_n=\int_{\R^N}e^{-V_n}dx.
\end{align*}
Then $\gamma_n$ becomes a log-concave probability measure. Therefore the diffusion process associated with the Dirichlet form $(\mathcal E_n, \mathcal F_n)$ in \eqref{defn: SDE} is a solution to the following SDE:
 \begin{align} \label{SDE: normal}
 dX^n_t=-\nabla V_n(X^n_t)dt+\sqrt{2}dW_t, \quad X_0=\x_n.
 \end{align}
 If
 $\gamma_n$ converges to a probability measure $\gamma_\infty$ weakly and $\x_n \to \x_\infty$, then $(\R^N, \|\cdot\|_2, \gamma_n, \x_n)$ converges to $(\R^N, \|\cdot\|_2, \gamma_\infty, \x_\infty)$ in the pmG sense and $(\R^N, \|\cdot\|_2, \gamma_\infty, \x_\infty)$ is an RCD$(0,\infty)$ space by the stability of the RCD property. Thus 
 the  \replaced{solution to SDE \eqref{SDE: normal}}{Brownian motions} on $H_n=(H,\|\cdot\|_n, \gamma_n, \x_n)$ converges weakly to 
 the diffusion associated with the Cheeger energy on the limit space (Theorem \ref{thm: mthm2}, \ref{thm: mthm2-0} and \cite[Theorem 1.5]{ASZ09}).
  \item[(b)]{\bf (SDE on \replaced{V}{v}ariable \replaced{C}{c}onvex \replaced{S}{s}ubsets with \replaced{V}{v}ariable \replaced{C}{c}onvex \replaced{P}{p}otentials)}
 Let $H=\R^N$ with $1<N<\infty$ and $U_n \subset \R^N$ be a smooth convex open set. We consider a convex functional $V_n \in C^{1,1}(U_n)$ and $V_n\equiv +\infty$ on $\R^N \setminus U_n$ with $\int_{U_n}e^{-V_n}dx<\infty$ for $n \in \N$.
 Take
 \begin{align*}
\gamma_n=\frac{1}{Z_n}e^{-V_n}dx|_{U_n}, \quad \text{whereby}\  Z_n=\int_{U_n}e^{-V_n}dx.
\end{align*}
 Then $\gamma_n$ becomes a log-concave probability measure. Therefore the diffusion process associated with the Dirichlet form $(\mathcal E_n, \mathcal F_n)$ in \eqref{defn: SDE} is a solution of the following SDE with reflection at the boundary:
 \begin{align} \label{sol: SDEB}
 dX^n_t=-\nabla V_n(X^n_t)dt+\sqrt{2}dW_t + \mathbf{n}_n(X_t)dL^n_t, \quad X_0=\x_n.
 \end{align}
Here $\mathbf{n}_n$ is an inner normal vector to $\partial U_n$ and $L^n$ is a continuous monotone non-decreasing process which increases only when $X_t \in \partial U_n$. 
 
 If $\gamma_n$ converges to a probability measure $\gamma_\infty$ weakly and $\x_n \to \x_\infty$, then $(\R^N, \|\cdot\|_2, \gamma_n, \x_n)$ converges to $(\R^N, \|\cdot\|_2, \gamma_\infty, \x_\infty)$ in the pmG sense and $(\R^N, \|\cdot\|_2, \gamma_\infty, \x_\infty)$ is an RCD$(0,\infty)$ space by the stability of the RCD property. 
 Thus 
 the \replaced{solution to}{diffusion processes associated with} \eqref{sol: SDEB} on $\overline{U}_n$  converges weakly to the diffusion associated with the Cheeger energy on the limit space (Theorem \ref{thm: mthm2}, \ref{thm: mthm2-0} and \cite[Theorem 1.5]{ASZ09}).
 \end{enumerate}
\end{exa}

 \section{Proof of Theorem \ref{thm: mthm2}} \label{sec: proof2}
 We first show the implication of (ii) $\implies$ (i) in Theorem \ref{thm: mthm2}.
 \\
 {\it Proof of (ii) $\implies$ (i) in Theorem \ref{thm: mthm2}}.
If we assume (ii), then it is obvious that the initial distributions $\m_n$ converge weakly to $\m_\infty$. Since the weak convergence of $\m_n$ to $\m_\infty$ is equivalent to the convergence of $m_n$ to $m_\infty$ in the sense of \eqref{eq: VGC} (easy to check), we finish the proof of the implication (ii) $\implies$ (i) in Theorem \ref{thm: mthm2}. \qed
\\
We now show the implication (i) $\implies$ (ii). 
\\
{\it Proof of (i) $\implies$ (ii) in Theorem \ref{thm: mthm2}}.
By Definition \ref{prop: Dconv}, there exist a complete separable metric space $(X,d)$ and a family of isometric embeddings $\iota_n: X_n \to X$ such that, for any bounded continuous function $f: X \to \R$ with bounded support, we have 
 \begin{align*} 
\int_{X}f d({\iota_n}_\#m_n) \to \int_{X}f d({\iota_\infty}_\#m_\infty).
 \end{align*}
Set the notation for the laws of Brownian motions as follows:
 $$\mathbb B^{\m_n}_n:=(\iota_n (B^n),  \mathbb P_n^{\m_n}), \quad \mathbb B^{\x_n}_n:=(\iota_n (B^n), \mathbb P_n^{\x_n}).$$ 
Hereafter we identify $\iota_n(X_n)$ with $X_n$ and we omit $\iota_n$ for simplifying the notation. 

To show the weak convergence of the Brownian motions, we have two steps. The first \added{step} is to show the weak convergence of finite-dimensional distributions, and the second is to show tightness. 
We first show the weak convergence of finite-dimensional distributions in the case that the initial distribution is the Dirac measure $\delta_{\x_n}$.
 \begin{lem} \label{lem: FDC1}{\bf (Convergence of Finite-Dimensional Distributions)}
For any $k \in \N$, $0 =t_0 < t_1 <t_2< \cdot \cdot \cdot <t_k<\infty$ and $f_1, f_2, ..., f_k \in C_b(X)$, the following holds:
\begin{align*} 
\mathbb E^{\x_n}[f_1(B^n_{t_1})\cdot\cdot\cdot f_k(B^n_{t_k})] \overset{n \to \infty}\to \mathbb E^{\x_\infty}[f_1(B^\infty_{t_1})\cdot\cdot\cdot f_k(B^\infty_{t_k})].
\end{align*} 
\end{lem}
\proof
Since\deleted{the} the limit Brownian motion $\mathbb B_\infty^{\x_\infty}$ is conservative,  it suffices to show the statement only for $f_1, f_2,...,f_k \in C_b(X) \cap L^2(X;m_\infty)$. In fact, for any $\e>0$ and $T>0$, there exists $R=R(\e,T)$ so that the open ball $B_R(\x_\infty)$ satisfies
$$\mathbb E^{\x_\infty}\1_{B_R(\x_\infty)}(B_t^\infty)=\mathbb P^{\x_\infty}(B_t^\infty \in B_R(\x_\infty)) \ge 1-\e, \quad \forall t \in [0,T].$$
If we know that $\mathbb E^{\x_n}(f(B_t^n))$ converges to $\mathbb E^{\x_\infty}(f(B_t^\infty))$ for any $ f\in C_b(X)  \cap L^2(X;m_\infty)$,  
then we know that 
\begin{align*}
 \lim_{n \to \infty}\mathbb P^{\x_n}(B_t^n \in B_R(\x_\infty)) =\lim_{n \to \infty} \mathbb E^{\x_n}(\1_{B_R(\x_\infty)}(B_t^n)) =\mathbb E^{\x_\infty}(\1_{B_R(\x_\infty)}(B_t^\infty)) \ge 1-\e,
\end{align*}
for any $t \in [0,T].$
Therefore, for any $f_1,...,f_k \in C_b(X)$, \added{and any small $\delta>0$, we can choose $R>0$ large enough so that} \deleted{we have}
\begin{align*}
&\lim_{n \to \infty} \mathbb E^{\x_n}(f_1(B_{t_1}^n)\cdot\cdot\cdot f_k(B_{t_k}^n))
\\
&= \lim_{n \to \infty} \mathbb E^{\x_n}\Bigl(f_1(B_{t_1}^n)\cdot\cdot\cdot f_k(B_{t_k}^n): \bigcap_{j=1}^k\{B_{t_j}^n \in B_R(\x_\infty)\}\Bigr) 
\\
& \qquad \quad +  \lim_{n \to \infty} \mathbb E^{\x_n}\Bigl(f_1(B_{t_1}^n)\cdot\cdot\cdot f_k(B_{t_k}^n): \Bigl(\bigcap_{j=1}^k\{B_{t_j}^n \in B_R(\x_\infty)\}\Bigr)^c\Bigr)
\\
&= \lim_{n \to \infty} \mathbb E^{\x_n}\Bigl(f_1\1_{B_R}(B_{t_1}^n)\cdot\cdot\cdot f_k\1_{B_R}(B_{t_k}^n)\Bigr) +{\delta}.
\end{align*}
Thus we may show the proof only for  $f_1, f_2,...,f_k \in C_b(X) \cap L^2(X;m_\infty)$.

Recall that we have the following equality (see Subsection \ref{subse: BM1}): for every $f \in C_b(X) \cap L^2(X;m_\infty)$,
\begin{align} \label{eq: KKT}
\mathbb E_n^{x}(f(B_t^n))=P^n_tf(x),
\end{align}
for {\it every} $x \in X_n$.
Here recall that $\{P^n_t\}_{t \ge 0}$ is the semigroup defined in \eqref{eq: HSBM} by the action of the heat flow whereby $P_t$ is a modification of the heat semigroup $H_t$ and $P^n_tf(x)$ can be defined for every point $x \in X^n$ if $f \in C_b(X) \cap L^2(X;m_\infty).$ Since the Brownian motion $(\{\mathbb P^x_n\}_{x \in X_n}, \{B_t^n\}_{t \ge 0})$ is constructed by the Kolmogorov extension theorem with the integral kernel $p_n(t,x,dy)$ of $\{P^n_t\}_{t \ge 0}$ \added{as in Section \ref{subse: BM1}}, the equality \eqref{eq: KKT} holds for every point $x \in X_n$

By using the Markov property, for all $n \in \EN$, we have 
\begin{align*} 
&\mathbb E_n^{\x_n}[f_{1}(B^{n}_{t_1})\cdot\cdot\cdot f_k(B^{n}_{t_k})]  \notag
\\
&=P^n_{t_1-t_0}\Bigl(f_{1} P^n_{t_2-t_1} \Bigl(f_{2}\cdot \cdot \cdot P^n_{t_k-t_{k-1}}f_{k} \Bigr)\Bigr)(\x_n) \notag
\\
&=:\mathcal P_k^n(\x_n).
\end{align*}
By \cite[Theorem 7.3]{AGMR15}, $\mathcal P_k^n$ is bounded Lipschitz on $X_n$ whose Lipschitz constant depends only on the curvature lower-bound $K$.

For later arguments, we extend $\mathcal P_k^n$ to the whole space $X$ (note that $\mathcal P_k^n$ is defined only on each $X_n$). The key point is to extend $\mathcal P_k^n$ to the whole space $X$ preserving its Lipschitz regularity and bounds.
\begin{prop}{\bf (\cite[Corollary 1,2]{Mc34})} \label{prop: McS}
Let $\widetilde{\mathcal P}_k^n$ be the function defined on the whole space $X$ as follows:
\begin{align} \label{eq: EX}
\widetilde{\mathcal P_k^n}(x) :=\Bigl( \sup_{a \in X_n}\{\mathcal P_k^n(a)-H{d(a,x)}\} \wedge \sup_{a \in X_n} \mathcal P_k^n(a) \Bigr)  \vee  \inf_{a \in X_n} \mathcal P_k^n(a), \quad x \in X.
\end{align}
Here $H$ denotes the same Lipschitz constant of the original function $\mathcal P_k^n$. 
Then $\widetilde{\mathcal P_k^n}$ is a bounded Lipschitz continuous function on the whole space $X$ with the same Lipschitz constant $H$ and the same bound. Moreover $\tilde{\mathcal P_k^n}=\mathcal P_k^n$ on the original domain $X_n$. The function $\widetilde{\mathcal P}_k^n$ is called {\it McShane extension of $\mathcal P_k^n$}. 
\end{prop}

We now return to the proof of Lemma \ref{lem: FDC1}. We have that 
\begin{align}
&\Bigl| \mathbb E_n^{\x_n}[f_{1}(B^{n}_{t_1})\cdot\cdot\cdot f_k(B^{n}_{t_k})] -\mathbb E_n^{\x_\infty}[f_{1}(B^{\infty}_{t_1})\cdot\cdot\cdot f_k(B^{\infty}_{t_k})]\Bigr| \notag
\\
&=|\mathcal P_k^n(\x_n)-{\mathcal P}_k^\infty(\x_\infty)| \notag
\\
& \le |{\mathcal P}_k^n(\x_n)-\tilde{\mathcal P_k^n}(\x_\infty)|+|\tilde{\mathcal P_k^n}(\x_\infty)-\mathcal P_k^\infty(\x_\infty)| \notag
\\
&=:({\rm I})_n+({\rm II})_n. \notag
\end{align}
Therefore it suffices to show (I)$_n \to 0$ and (II)$_n \to 0$ as $n \to \infty$. 

We first discuss to show (I)$_n \to 0$.  Since 
$\|P_t^nf\|_\infty =\|f\|_\infty \|\int_{X_n}p_n(t,x,y)m_n(dy)\|_\infty \le \|f\|_\infty,$ for any $f \in C_b(X_n) \cap L^2(X;m_\infty)$, we have 
\begin{align} \label{ineq: BDDSG}
\sup_{n \in \N}\|\mathcal P_k^n\|_\infty \le \prod_{i=1}^k \|f_i\|_\infty<\infty.
\end{align}
Therefore, by Proposition \ref{prop: McS}, it holds that 
\begin{align} \label{uniform bound of HS}
\sup_{n \in \N}\|\tilde{\mathcal P_k^n}\|_\infty<\infty.
\end{align}
By \cite[Theorem 7.3]{AGMR15},  we have that ${\rm Lip}_X(P_t^nf) \le C(t, K) \|f\|_\infty$ for any $f \in L^\infty(X_n;m_n) \cap L^2(X;m_\infty)$ for some positive $C(t,K)$ depending only on $t, K$. Here ${\rm Lip}_X(f)$ means the global Lipschitz constant of a Lipschitz function $f$ on $X$. Thus by considering \eqref{ineq: BDDSG}, there exists a constant $L$ depending only on $t_k, K$ and $\|f_1\|_\infty,....,\|f_k\|_\infty$ (but independent of $n$) so that 
 \begin{align*}
  \sup_{n \in \N}{\rm Lip}_X(\mathcal P_k^n)  &\le  \sup_{n \in \N}C(t_k,K)\| f_k\mathcal P_{k-1}^n\|_\infty <L<\infty. \notag
  \end{align*} 
By the property of the McShane extension in Proposition \ref{prop: McS}, we have that 
   \begin{align} \label{ineq: MCSG}
  \sup_{n \in \N}{\rm Lip}_X(\tilde{\mathcal P_k^n})  &\le  \sup_{n \in \N}C(t_k,K)\| f_k\mathcal P_{k-1}^n\|_\infty <L<\infty. 
  \end{align} 
Thus we have 
\begin{align*}
({\rm I})_n= |{\mathcal P_k^n}(\x_n)-\tilde{\mathcal P_k^n}(\x_\infty)|
&= |\tilde{\mathcal P_k^n}(\x_n)-\tilde{\mathcal P_k^n}(\x_\infty)|
\\
&\le {\rm Lip}(\tilde{\mathcal P_k^n})d(\x_n,\x_\infty)
\\
&\le Ld(\x_n,\x_\infty)
\\
& \to 0 \quad (n \to \infty).
\end{align*}

We next show that (II)$_n \to 0$. By \eqref{uniform bound of HS} and \eqref{ineq: MCSG}, we can apply the Ascoli--Arzel\'a theorem to $\{\tilde{\mathcal P_k^n}\}_{n \in \N}$ so that $\{\tilde{\mathcal P_k^n}\}_{n \in \N}$  is relatively compact. Therefore, for any subsequence $\{\tilde{\mathcal P_k^{n'}}\}_{\{n'\}}$ whereby $\{n'\} \subset \{n\}$, there exists a further subsequence $\{\tilde{\mathcal P_k^{n''}}\}_{\{n''\}}$ whereby $\{n''\} \subset \{n'\}$ satisfying 
\begin{align} \label{convergence: HS}
\tilde{\mathcal P_k^{n''}} \to F'' \quad \text{uniformly in}\ X.
\end{align}

On the other hand, we have that ${\mathcal P}_k^n$ converges to ${\mathcal P}_k^\infty$ $L^2$-strongly in the sense of Definition \ref{defn: Weak}. 
We give a proof below. 
\begin{lem} \label{lem: SCCHS}
${\mathcal P}_k^n$ converges to ${\mathcal P}_k^\infty$ in the $L^2$-strong sense in Definition \ref{defn: Weak}. 
\end{lem}
\proof
By Theorem \ref{thm: Mosco of SG}, the statement is true for $k=1$. Assume that the statement is true when $k=l$. Since we have 
$$\mathcal P_{l+1}^n=P^n_{t_{l+1}-t_l}(f_{l+1}^{(n)}\mathcal P_l^n),$$
by Theorem \ref{thm: Mosco of SG}, it is sufficient to show $f_{l+1}\mathcal P_l^n \to f_{l+1}\mathcal P_l^\infty$ strongly in $L^2$. This is obvious to be true because $\mathcal P_l^n \to \mathcal P_l^\infty$ strongly (the assumption of the induction), $f_{l+1} \in C_b(X)$ and $\mathcal P_l^n$  is bounded uniformly in $n$ thanks to \eqref{uniform bound of HS}. Thus the statement is true for any $k \in \N$. 
\qed

We return to the proof of Lemma \ref{lem: FDC1}.
\\
{\it Proof of Lemma \ref{lem: FDC1}.}
By using Lemma \ref{lem: SCCHS} and \eqref{convergence: HS}, it is obvious to check that 
$$F''|_{X_\infty}=\mathcal P_k^\infty,$$
whereby $F''|_{X_\infty}$ means the restriction of $F''$ into $X_\infty$. 
The R.H.S. $\mathcal P_k^\infty$ of the above equality is clearly independent of choices of subsequences and thus the limit $F''|_{X_\infty}$ is independent of choices of subsequences. Thus we conclude that 
\begin{align} \label{concl: UC}
\tilde{\mathcal P_k^n} \to \mathcal P_k^\infty \quad \text{uniformly in } \ X_\infty.
\end{align}
Now we return to show (II)$_n$ goes to zero. By \eqref{concl: UC}, we have that 
\begin{align*}
({\rm II})_n=|\tilde{\mathcal P_k^n}(\x_\infty)-\mathcal P_k^\infty(\x_\infty)|
&\le \|\tilde{\mathcal P_k^n}-\mathcal P_k^\infty\|_{\infty, X_\infty}
\\
& \to 0 \quad (n \to \infty).
\end{align*}
Here $\|\cdot\|_{\infty, X_\infty}$ means the uniform norm on $X_\infty$. Thus we finish the proof of Lemma \ref{lem: FDC1}.
\qed

We next show the weak convergence of finite-dimensional distributions for the case that initial distributions are $W_1$-convergent, which includes $\m_n$ for the case of $m_n(X_n)=\infty$. 
\begin{lem} \label{lem: FDC1-2}
Let $\{\nu_n\}_{n \in \N} \subset \mathcal P(X_n)$ be a sequence of probability measures on $X_n \subset X$ converging to $\nu_\infty \in \mathcal P(X_\infty)$ in $W_1$-distance.
Then, for any $k \in \N$, $0 =t_0 < t_1 <t_2< \cdot \cdot \cdot <t_k<\infty$ and $f_1, f_2, ..., f_k \in C_b(X) \cap L^2(X;m_\infty)$, the following holds:
\begin{align*}
\mathbb E^{\nu_n}[f_1(B^n_{t_1})\cdot\cdot\cdot f_k(B^n_{t_k})] \overset{n \to \infty}\to \mathbb E^{\nu_\infty}[f_1(B^\infty_{t_1})\cdot\cdot\cdot f_k(B^\infty_{t_k})].
\end{align*} 
\end{lem}
\proof
By the same argument at the beginning of Lemma \ref{lem: FDC1}, it suffices to show the statement for any $f_1, f_2,...,f_k \in C_b(X) \cap L^2(X;m_\infty)$.
Recall that we set in Lemma \ref{lem: FDC1} as follows:
\begin{align*} 
&\mathbb E_n^{x}[f_{1}(B^{n}_{t_1})\cdot\cdot\cdot f_k(B^{n}_{t_k})]  \notag
\\
&=P^n_{t_1-t_0}\Bigl(f_{1}^{(n)} P^n_{t_2-t_1} \Bigl(f_{2}^{(n)} \cdot \cdot \cdot P^n_{t_k-t_{k-1}}f_{k}^{(n)} \Bigr)\Bigr)(x) \notag
\\
&=:\mathcal P_k^n(x).
\end{align*}
By the Kantorovich--Rubinstein duality (see e.g., \cite[Theorem 5.10]{V09}), we have 
\begin{align*}
W_1(\nu_n,\nu_\infty)=\frac{1}{L}\sup\{\int_{X}f d\nu_n-\int_X f\nu_\infty: f \in {\rm Lip}_b(X), \ {\rm Lip}_X(f) \le L\}.
\end{align*}
According to \eqref{uniform bound of HS} and \eqref{ineq: MCSG}, we have that $\tilde{\mathcal P_k^n}$ is bounded and $\sup_{n \in \N}{\rm Lip}(\tilde{\mathcal P_k^n})<L<\infty$ for some constant $L$. Thus we have that 
\begin{align} \label{duality: KR}
 \Bigl| \int_{X}\tilde{\mathcal P}_k^n\ d\nu_n - \int_{X}\tilde{\mathcal P}_k^n d\nu_\infty\Bigr| \le LW_1(\nu_n,\nu_\infty). 
\end{align}
Since $\tilde{\mathcal P_k^n}$ converges to $\mathcal P_k^\infty$ uniformly in $C_b(X_\infty)$ by \eqref{concl: UC}, and $\nu_n$ converges to $\nu_\infty$ in the $W_1$-distance, by using \eqref{duality: KR}, we have that 
\begin{align*}
& \Bigl| \mathbb E^{\nu_n}[f_1(B^n_{t_1})\cdot\cdot\cdot f_k(B^n_{t_k})]  - \mathbb E^{\nu_\infty}[f_1(B^\infty_{t_1})\cdot\cdot\cdot f_k(B^\infty_{t_k})] \Bigr|
\\
&= \Bigl| \int_{X}\mathcal P_k^n\ d\nu_n - \int_{X}\mathcal P_k^\infty\ d\nu_\infty \Bigr|
\\
&\le  \Bigl| \int_{X}\mathcal P_k^n\ d\nu_n - \int_{X}\tilde{\mathcal P_k^n} d\nu_\infty\Bigr|+ \Bigr| \int_{X}\tilde{\mathcal P_k^n} d\nu_\infty- \int_{X}\mathcal P_k^\infty\ d\nu_\infty \Bigr|
\\
&\le  LW_1(\nu_n,\nu_\infty)+ \|\tilde{\mathcal P_k^n}-\mathcal P_k^n\|_{\infty, X_\infty}\int_{X}\ d\nu_\infty
\\
& \to 0, \quad n \to \infty.
\end{align*}
Thus we have completed the proof.
\qed

We now show the weak convergence of finite-dimensional distributions for the case that initial distributions are $\frac{1}{m_n(X_n)}m_n$, which corresponds to the case of $m_n(X_n)<\infty$.
\begin{lem} \label{lem: FDC-2}
Let $m_n(X_n)<\infty$ for any $n\in \N$.
Then, for any $k \in \N$, $0 =t_0 < t_1 <t_2< \cdot \cdot \cdot <t_k<\infty$ and $f_1, f_2, ..., f_k \in C_b(X)$, the following holds:
\begin{align*} 
\mathbb E^{\m_n}[f_1(B^n_{t_1})\cdot\cdot\cdot f_k(B^n_{t_k})] \overset{n \to \infty}\to \mathbb E^{\m_\infty}[f_1(B^\infty_{t_1})\cdot\cdot\cdot f_k(B^\infty_{t_\infty})].
\end{align*} 
\end{lem}
\proof
Because of $m_n(X_n)<\infty$, we have $f \in L^2(X, m_n)$ for all $f \in C_b(X)$ for any $n \in \EN$. 
Since $\m_n$ converges weakly to $\m_\infty$ in $\mathcal P(X)$, for any $\e>0$, there exists a compact set $ K \subset X$ so that 
$$\sup_{n \in \N}\m_n(K^c)<\e.$$
Thus, by \eqref{ineq: BDDSG}, for any $\delta>0$, there exists a compact set $ K \subset X$ so that 
\begin{align} \label{ineq: SFDC}
\sup_{n \in \N}\Bigl| \int_{X_n}\mathcal P_k^n d\m_n-\int_{K}\mathcal P_k^n\ d\m_n\Bigr| \le  \Bigl(\prod_{i=1}^k\|f_i\|_\infty\Bigr) \sup_{n \in \N}\m_n(K^c)< \delta.
\end{align}
Take $r>0$ so that $K \subset B_r(\x_n):=\{x \in X: d(\x_n,x)<r\}$. Let $\tilde{\1}_r^R$ denote the following function: ($r<R$)
\begin{align*}
\tilde{\1}_r^R(x)=
\begin{cases} \dis
1, &\text{$x \in B_r(\x_n)$},
\\
\dis 1-\frac{d(x, B_r(\x_n))}{R-r}, &\text{$x \in B_R(\x_n) \setminus B_r(\x_n)$},
\\
0, &\text{otherwise.}
\end{cases}
\end{align*}
Then $\tilde{\1}_r^R \in C_{bs}(X)$.
Thus, by Theorem \ref{thm: Mosco of SG} and \eqref{ineq: SFDC}, for any $\delta>0$, there exists $r>0$ so that
\begin{align*}
& \Bigl| \mathbb E^{\m_n}[f_1(B^n_{t_1})\cdot\cdot\cdot f_k(B^n_{t_k})]  - \mathbb E^{\m_\infty}[f_1(B^\infty_{t_1})\cdot\cdot\cdot f_k(B^\infty_{t_\infty})] \Bigr|
\\
&=  \Bigl| \int_{X_n}\mathcal P_k^n d\m_n - \int_{X_\infty}\mathcal P_k^\infty d\m_\infty \Bigr|
\\
&=  \Bigl| \int_{X_n}\mathcal P_k^n d\m_n - \int_{X_n}\tilde{\1}_r^R \mathcal P_k^n d\m_n+ \int_{X_n}\tilde{\1}_r^R\mathcal P_k^n d\m_n
\\
&\quad - \int_{X_n}\tilde{\1}_r^R\mathcal P_k^\infty d\m_\infty
+\int_{X_n} \tilde{\1}_r^R \mathcal P_k^\infty d\m_\infty- \int_{X_\infty}\mathcal P_k^\infty d\m_\infty \Bigr|
\\
& \le \delta+ \Bigl|\int_{X}\tilde{\1}_r^R \mathcal P_k^n d\m_n- \int_{X}\tilde{\1}_r^R  \mathcal P_k^\infty d\m_\infty\Bigr| + \delta
\\
& \overset{n \to \infty}\to 2\delta. 
\end{align*} 
In the fourth line, the first $\delta$ comes from using \eqref{ineq: SFDC} and the second $\delta$ comes from using the tightness of the single measure $m_\infty$. The the middle term in the fourth line converges to zero thanks to the $L^2$-strong convergence of the heat semigroup $P_t$ in the sense of Definition \ref{defn: Weak}. 
Note that the total mass $m_n(X_n) \to m_\infty(X_\infty)$\added{(}$\le \infty$\added{)} because of the pmG convergence. 
Thus we have completed the proof.
\qed

Now we show the tightness of $\{\mathbb B^{\m_n}\}$. For later arguments, we show the tightness for more general initial distributions $\nu_n$ than $\m_n$.
\begin{lem} \label{lem: T1}
Let $\nu_n \in \mathcal P(X_n)$ satisfy the following conditions:
\begin{enumerate}
\item[(i)] $\nu_n \to \nu_\infty$ weakly in $\mathcal P(X)$;
\item[(ii)] $\nu_n$ is absolutely continuous with respect to $m_n$ with $d\nu_n=\phi_ndm_n$ and there exists a positive constant $M$ so that, for any $r>0$, 
$$\sup_{n \in \N}\|\phi_n\|_{\infty, B_r(\x_n)} <M<\infty.$$
\end{enumerate}
Then 
 $\{\B^{\nu_n}\}_{n \in \N}$ is tight in $\mathcal P(C([0,\infty), X))$.
 \end{lem}
 \proof
Let us denote the law of $h(B^n)$ for $h \in {\rm Lip}_b(X)$ as follows:
$$
 \mathbb B^{\nu_n, h}=(h(B^n), \mathbb P_n^{\nu_n}).
$$
It is easy to show that ${\rm Lip}_b(X)$ strongly separates points in $C_b(X)$, that is, for every $x$ and $\e>0$, there exists a finite set $\{h_i\}_{i=1}^l \subset {\rm Lip}_b(X)$ so that 
$$\inf_{y: d(y,x) \ge \e}\max_{1 \le i \le l}|h_i(x)-h_i(y)|>0.$$
Therefore, by \cite[Corollary  3.9.2]{EK86} with Lemma \ref{lem: FDC1-2}, the following two statements are equivalent:
\begin{enumerate}
\item[(i)]  $\{\B^{\nu_n}\}_{n \in \N}$ is tight in $\mathcal P(C([0,\infty), X))$;
\item[(ii)]  $\{\B^{\nu_n,h}\}_{n \in \N}$ is tight in $\mathcal P(C([0,\infty), \R))$.
\end{enumerate}
Thus we will show that, for any $h \in {\rm Lip}_b(X)$, 
 \begin{align*}  
\{\B^{\nu_n,h}\}_{n \in \N} \quad \text{ is tight in } \quad \mathcal P(C([0,\infty); \R)).
 \end{align*}
We note that, although \cite[Corollary 3.9.2]{EK86} gives sufficient conditions for tightness only in the c\`adl\`ag space $D([0,\infty);X)$, since the law\deleted{s} of each Brownian motion\deleted{s} $\mathbb B_n^{\m_n}$ \added{for $n \in \EN$} \replaced{has}{have} \replaced{its}{their} support on the space of continuous paths $C([0,\infty);X)$, the tightness in $D([0,\infty);X)$ implies the tightness in $C([0,\infty), X)$. See, e.g., \cite[Lemma 5 in Appendix]{FK97} for this point.

Since $\nu_n$ converges weakly to $\nu_\infty$ in $\mathcal P(X)$, 
\added{the set of} the laws of the initial distributions $\{{(h(B^n_0)}, \mathbb P_n^{\nu_n})\}_{n \in \N}=\{h_\#\m_n\}_{n \in \N}$ is clearly tight in $\mathcal P(\R)$. For $\delta>0$, let us define
 $$
 L_{\eta,T}^{n,h}(x):=\mathbb P_n^x(\sup_{\substack{0 \le s,t\le T \\ |t-s|\le \eta}}|h(B_t)-h(B_s)| >\delta).
 $$
The desired result we would like to show is the following:
 \begin{align} \label{convergence: tightness-RCD}
 \lim_{\eta\to 0}\sup_{n \in \N}\int_{X_n}L_{\eta,T}^{n,h}\ d\nu_n=0,
 \end{align}
 for any $T>0$.
By conditions (i) and (ii) in this lemma, for any $\e>0$, there exists $R>0$ so that 
 \begin{align*}
 \int_{X_n}L_{\eta,T}^{n,h}d\nu_n
 &=\|\phi_n\1_{B_R(\x_n)}\|_\infty \int_{X_n}L_{\eta,T}^{n,h}\1_{B_R(\x_n)} dm_n+\nu_n(B_R^c(\x_n))
 \\
 &<M\int_{X_n}L_{\eta,T}^{n,h}\1_{B_R(\x_n)} dm_n+\e.
 \end{align*}
It suffices to show, for any $T,R>0$, 
$$\lim_{\eta \to 0}\sup_{n \in \N}\int_{X_n} L_{\eta,T}^{n,h}\1_{B_R(\x_n)} dm_n=0.$$
Let $m_{n,R}:=\1_{Y^R_n}m_n$ whereby 
$$Y^R_n=\overline{B_R(\x_n)}$$
 is the closure of the open ball $B_R(\x_n)$.
We have
\begin{align*} 
\int_{X_n}L_{\eta,T}^{n,h}\ dm_{n,R}  & =\mathbb P_{n, R+r}^{m_{n,R}}\Bigl(\sup_{\substack{0 \le s,t\le T \\ |t-s|\le \eta}}|h(B^n_t)-h(B^n_s)| >\delta: \Lambda_r \Bigr)
\\
&\quad +\mathbb P^{m_{n,R}}\Bigl(\sup_{\substack{0 \le s,t\le T \\ |t-s|\le \eta}}|h(B^n_t)-h(B^n_s)| >\delta: \Lambda_r^c \Bigr) \notag
\\
&:= {\rm (I)_{n,\eta}} + {\rm (II)_{n,\eta}},
\end{align*}
whereby $\Lambda_r:=\{w \in \Omega^n: \sup_{0 \le t \le T}|d_n(B^n_t,\x_n)-d_n(B^n_0,\x_n)|<r\}$. Here $\mathbb P_{n, r}^x$ is a conservative diffusion process associated with $(\C_n^{r}, \mathcal F_n^{r})$
$${\sf Ch}_n^r(f)=\frac{1}{2}\int_{Y^r_n}|\nabla f|^2_{w,Y^r_n} dm_{n,r}, \quad \mathcal F_n^r:=\{f \in L^2(Y^r_n;m_{n,r}): {\sf Ch}^r_n(f)<\infty \}.$$
Recall that $|\nabla f|^2_{w,Y^r_n}$ means the minimal weak upper gradient on $Y^r_n$ (see Subsection \ref{subsub: Ch}). 
We note that the Cheeger energy $\C_n^{r}$ on the closed ball $Y^r_n$ is also quadratic because of \cite[Theorem 4.19]{AGS14b}.
Since closed balls are not necessarily convex subset in $X_n$, the closed ball $Y^r_n$ is not necessarily an RCD$(K,\infty)$ space. However, we can still construct the Brownian motion on $Y^r_n$ since we have that $(\C_n^{r}, \mathcal F_n^{r})$ is quadratic (\cite[Theorem 4.19]{AGS14b}) and $[d(x, \cdot)] \le m_{n,r}$ (\cite[(iv) Theorem 4.18]{AGS14b}) for any fixed $x \in Y^r_n$, which imply that $(\C_n^{r}, \mathcal F_n^{r})$ becomes a quasi-regular Dirichlet form by the same manner of \cite[Lemma 6.7]{AGS14b} and \cite[Theorem 1.2]{ASZ09} (see also \cite[\S 7.2]{AGMR15}). Here $[f]$ means the energy measure of the Cheeger energy (see \cite[(4.21)]{AGS14b}) and $[d(\x_n, \cdot)] \le m_{n,r}$ means 
$$\frac{d[d(\x_n,\cdot)]}{dm_{n,r}}(y) \le 1 \quad \text{$m_{n,r}$-a.e. $y \in Y^r_n$}.$$ 
Note that although \cite[Lemma 6.7]{AGS14b} assumed the RCD$(K,\infty)$ condition, only the quadraticity of the Cheeger energy and $[d(\x_n, \cdot)] \le m_{n,r}$ are used  to construct the Brownian motions, and the CD$(K,\infty)$ condition is not necessary (see also \cite[\S 4]{LS17} for more detailed studies of the Cheeger energies and Brownian motions on subsets in RCD$(K,\infty)$ spaces).

We first estimate (I)$_{n,\eta}$.
By {Lyons-Zheng decomposition} (\cite{LZ94}, and see also \cite[Section 5.7]{FOT11}), we have 
\begin{align*} 
h(B^n_t)-h(B^n_s)=\frac{1}{2}(M_t^{[h]}-M_s^{[h]})+\frac{1}{2}(M_{T-t}^{[h]}(r_T)-M_{T-s}^{[h]}(r_T)), \quad \text{$\mathbb P_{R+r}^{m_{n,R+r}}$-a.e.},
\end{align*}
for $ 0 \le t \le T.$ 

Then by time-symmetry (see \cite[Lemma 5.7.1]{FOT11}), we have 
\begin{align} \label{ineq: tight3}
{\rm (I)}_{n,\eta}& \le \mathbb P_{R+r}^{m_{n,R+r}}(\sup_{\substack{0 \le s,t\le T \\ |t-s|\le \eta}}|h(B^n_t)-h(B^n_s)| >\delta) \notag
\\
&\le \mathbb P_{R+r}^{m_{n,R+r}}(\sup_{\substack{0 \le s,t\le T \\ |t-s|\le \eta}}\bigl| M_t^{[h],n}-M_s^{[h],n} \bigr| > \delta) \notag
\\
& \quad + \mathbb P_{R+r}^{m_{n,R+r}}(\sup_{\substack{0 \le s,t\le T \\ |t-s|\le \eta}}\bigl| M_{T-t}^{[h],n}(r_T)-M_{T-s}^{[h],n}(r_T) \bigr| > \delta) \notag
 \\
& = 2\mathbb P_{R+r}^{m_{n,R+r}}(\sup_{\substack{0 \le s,t\le T \\ |t-s|\le \eta}} \bigr| M_t^{[h],n}-M_s^{[h],n} \bigr|  > \delta). 
\end{align}
Since $M^{[h],n}$ is a continuous martingale, by the martingale representation theorem, there exists the one-dimensional Brownian motion  ${\mathbf B}^n(t)$ on an extended probability space $(\tilde{\Omega}, \tilde{\mathcal M}, \tilde{\mathbb P}_n^x)$ whereby $M^{[h],n}$ is represented as a time-changed Brownian motion with respect to the quadratic variation $\tilde{\mathbb P}_n^x$-a.s, q.e.\ $x \in Y^{R+r}_n$ (see, e.g.,  Ikeda--Watanabe \cite[Chapter II Theorem 7.3']{IW89}). That is, for q.e.\ $x \in Y^{R+r}_n$,
\begin{align*}
M^{[h],n}_t={\mathbf B}^n(\la M^{[h],n} \ra_t)={\mathbf B}^n\Bigl(\int_0^t \frac{d\mu^n_{\la h \ra}}{dm_n}(B_u^n)du\Bigr)={\mathbf B}^n\Bigl(\int_0^t |\nabla h|^2_{w, Y^{R+r}_n}(B_u^n)du\Bigr)  \quad \text{$\tilde{\mathbb P}_n^x$-a.s.}
\end{align*}
The last equality followed from \cite[(iv) Theorem 4.18]{AGS14b}.
Since $|\nabla h|_{w, Y^{R+r}_n} \le {\rm Lip}(h)$, we have
\begin{align*}
&\{\omega \in \tilde{\Omega}:\sup_{\substack{0 \le s,t\le T \\ |t-s|\le \eta}}\bigr| M_t^{[h],n}-M_s^{[h],n} \bigr|  > \delta \} \notag
\\
&=\{\omega \in \tilde{\Omega}:\sup_{\substack{0 \le s,t\le T \\ |t-s|\le \eta}}\Bigr| {\mathbf B}^n\Bigl(\int_0^t |\nabla h|^2_{w, Y^{R+r}_n}(B_u^n)du\Bigr)-{\mathbf B}^n\Bigl(\int_0^s |\nabla h|^2_{w, Y^{R+r}_n}(B_u^n)du\Bigr) \Bigr|  > \delta \} \notag
\\
&\subset \{\omega \in \tilde{\Omega}:\sup_{\substack{0 \le s,t\le {\rm Lip}(h)^2T \\ |t-s|\le {\rm Lip}(h)^2\eta}}\bigr| {\mathbf B}^n(t)-{\mathbf B}^n(s) \bigr|  > \delta \}. \notag
\end{align*}
Let $\mathbb W$ be the standard Wiener measure on $C([0,\infty);\R)$. Let 
$$\theta(\eta,h):=\mathbb W_n(\sup_{\substack{0 \le s,t\le {\rm Lip}(h)^2T \\ |t-s|\le {\rm Lip}(h)^2\eta}}|\omega(t)-\omega(s)|>\delta).$$
By \eqref{ineq: tight3} and noting $\sup_{n \in \N}m_n(B_{R+r}(\x_n))<\infty$ because of the weak convergence of $m_n$, we have, for any $T>0$,
\begin{align} \label{tight: est1}
 {\rm (I)}_{n, \eta}& \le \sup_{n \in \N}\int_{X_n}L_{\eta,T}^{n,h} dm_{n, R+r} \notag
 \\
 &\le \sup_{n \in \N}2\mathbb P_{R+r}^{m_{n,R+r}}(\sup_{\substack{0 \le s,t\le T \\ |t-s|\le \eta}} \bigr| M_t^{[h]}-M_s^{[h]} \bigr|  > \delta) \notag
 \\
 &\le 2\theta(\eta,h)\sup_{n \in \N}m_{n}(B_{R+r}(\x_n)) \notag
 \\
 & \overset{\eta \to 0}\to 0.
 \end{align}
 
We now estimate (II)$_{n,\eta}$. We have the following estimate:
\begin{align} \label{tight: est2}
{\rm (II)}_{n, \eta}&=\mathbb P^{m_{n,R}}\Bigl(\sup_{\substack{0 \le s,t\le T \\ |t-s|\le \eta}}|h(B^n_t)-h(B^n_s)| >\delta: \Lambda_r^c \Bigr) \notag
\\
& \le 6m_{n}(B_{R+r}(\x_n))\frac{1}{\sqrt{2\pi}}\int_{\frac{2r}{3\sqrt{{\rm Lip}(h)^2T}}}^\infty \exp\{-\frac{s^2}{2}\}ds \notag
\\
& \le c\exp\{{{\sf c_2}(R+r)^2}\}\int_{\frac{2r}{3\sqrt{{\rm Lip}(h)^2T}}}^\infty \exp\{-\frac{s^2}{2}\}ds. \notag
\\
& \le c\exp\{{\sf c_2}(R+r)^2\} \frac{3\sqrt{{\rm Lip}(h)^2T}}{2r}\exp\{-\frac{r^2}{18{\rm Lip}(h)^2T}\} \notag
\\
& \overset{r \to \infty}\to 0.
\end{align}
Here $c>0$ is a constant independent of $n$. 
In the second line above, we used \cite[Lemma 5.7.2]{FOT11}, in the third line, we used the volume growth estimate \eqref{asmp: basic1} and, in the \replaced{fourth}{forth} line, we used the fact $\int_x^\infty\exp\{\frac{s^2}{2}ds\} \le \frac{1}{x}\exp\{-\frac{x^2}{2}\}$. 
Thus, by \eqref{tight: est1} and \eqref{tight: est2}, we have that, for any $R>0$,
$$\lim_{\eta \to 0}\sup_{n \in \N}\int_{X_n} L_{\eta,T}^{n,h}\1_{B_R(\x_n)} dm_n=\lim_{\eta \to 0}\sup_{n \in \N}\Bigl({\rm (I)}_{n,\eta}+{\rm (II)}_{n,\eta}\Bigr)=0.$$
Thus we have the desired result \eqref{convergence: tightness-RCD}.
\qed
\\
We resume to prove Theorem \ref{thm: mthm2}.
\\
{\it Proof of Theorem \ref{thm: mthm2}}. 
 It is easy to check that conditions (i) and (ii) in Lemma \ref{lem: T1} are satisfied with $\nu_n=\m_n$ in the both cases of $m_n(X_n)=\infty$ and $m_n(X_n)<\infty$. Thus we have shown the tightness. 
By using Lemma \ref{lem: FDC-2}, we have completed the proof of (i) $\implies$ (ii) in Theorem \ref{thm: mthm2} in the case of $m_n(X_n)<\infty$.
Moreover, we can check easily that the conditions in Lemma \ref{lem: FDC1-2} are satisfied with $\nu_n=\m_n$ in the case of $m_n(X_n)=\infty$ (see \cite[Remark 4.6]{GMS13}). 
Therefore, we have completed the proof of (i) $\implies$ (ii) in Theorem \ref{thm: mthm2} in the case of $m_n(X_n)=\infty$.
We finish the proof of (i) $\implies$ (ii) in Theorem \ref{thm: mthm2}. \qed

 \section{Proof of Theorem \ref{thm: mthm2-0}} \label{subsec: I-III}
We show the following statement: for any $\e>0$, 
\begin{itemize}
 \item[{{\bf (iii)}$_{\ge \e}$}\deleted{\replaced{{\bf (iii)}'$_{>0}$}{{\bf (iii)}$_\e$}}]  
 There exist a complete separable metric space $(X,d)$ and isometric embeddings $\iota_n: X_n \to X\ (n \in \EN)$ so that
 \begin{align*}
(\iota_n(B^n), \mathbb P_n^{\x_n}) \ {\to} \ (\iota_\infty(B^\infty), \mathbb P_\infty^{\x_\infty}) \quad \text{weakly} \quad \text{in $\mathcal P(C([\e,\infty); X))$}.
\end{align*}
\end{itemize}
We first discuss the case of\deleted{the} condition {\bf (A)}, that is, $m_n(X_n)=1$.
\\
Since we have already shown the weak convergence of the finite-dimensional distributions under the general RCD$(K,\infty)$ condition for starting points $\x_n$ in Lemma \ref{lem: FDC1}, it suffices to prove the tightness:
\begin{lem} \label{lem: T1-1} \added{Under condition {\bf (A)},}
 $\{\B^{\x_n}_n\}_{n \in \N}$ is tight in $\mathcal P(C([\e,\infty), X))$ for any $\e>0$.
 \end{lem}
 \proof 
In the proof of \cite[Theorem 7.7]{GMS13}, we have 
\begin{align*} 
\sup_{n \in \N} {\rm Ent}_{m_n}(p_n({\e},\x_n,dy))=\sup_{n \in \N} {\rm Ent}_{m_n}(\mu_\e^{n, \x_n})<\infty,
\end{align*}
\added{where $\mu_\e^{n,\x_n}:=\mathcal H^n_\e\delta_{\x_n}$ defined in Subsection \ref{subsub: RE}.}
Let $\B^{\x_n}_n$ and $\B^{m_n}_n$ be restricted to the path space $C([\e,\infty), X)$.
By using Markov property, we have that 
$$\frac{d{\mathbb B_n^{\x_n}}}{d\mathbb B_n^{m_n}}=p(\e,\x_n, B_\e^n).$$
In fact, we have that, for any Borel measurable functions $F: C([\e,\infty), X): \to \R,$
\begin{align*}
\mathbb E^{\x_n}(F(B^n_{\e+\cdot}))&=\mathbb E^{\x_n}(\mathbb E^{B_\e^n}(F))
\\
&=\int_{X_n}\mathbb E^y\bigl(F(B^n_\cdot)\bigr)p_n(\e,\x_n, dy) 
\\
&= \int_{X_n}\mathbb E^y(F(B^n_\cdot))p_n(\e,\x_n,y)m_n(dy)
\\
&=\mathbb E^{m_n}(p_n(\e,\x_n, B^n_0)F(B^n_\cdot))
\\
&=\mathbb E^{m_n}(p_n(\e,\x_n, B^n_\e)F(B^n_{\e+\cdot})),
\end{align*}
whereby in the last line, we used the stationarity. 

Let us denote ${\rm Ent}_\nu(\mu)={\rm Ent}(\mu|\nu)$. By the fact that $\sup_{n \in \N}{\rm Ent}_{m_n}(p_n(\e,\x_n,dy))<\infty$, we have
\begin{align*}
 \sup_{n \in \N} {\rm Ent}(\mathbb B_n^{\x_n}|\mathbb B_n^{m_n})&= \sup_{n \in \N} \int_{\Omega}p_n(\e,\x_n, B_\e^n)\log\Bigl\{p_n(\e,\x_n, B_\e^n)\Bigr\}d\mathbb P_n^{m_n} 
 \\
 &= \sup_{n \in \N} \int_{X_n}P^n_\e\Bigl(p_n(\e,\x_n, \cdot)\log\{p_n(\e, \x_n, \cdot)\}\Bigr)dm_n
 \\
 &= \sup_{n \in \N}\int_{X_n}p_n(\e,\x_n, \cdot)\log\{p_n(\e, \x_n, \cdot)\}dm_n
 \\
 &= \sup_{n \in \N}{\rm Ent}_{m_n}(p_n(\e,\x_n,dy))<\infty.
 \end{align*}
In the third line above, we used the invariance property of $m_n$ with respect to the heat semigroup $\{P^n_t\}_{t \ge 0}$ whereby 
$$\int_{X_n}P^n_tf dm_n=\int_{X_n}fdm_n.$$
Since $\{\mathbb B_n^{m_n}\}_{n \in \N}$ is tight by Lemma \ref{lem: T1}, by using the tightness criterion with respect to the entropy \cite[Proposition 4.1]{GMS13}, we have the tightness of $\{\mathbb B_n^{\x_n}\}_{n \in \N}$. 
 \qed

{\it Proof of Theorem \ref{thm: mthm2-0} \added{in the case of {\bf (A)}}}.
By the weak convergence of the finite-dimensional distributions in Lemma \ref{lem: FDC1}, and the tightness in Lemma \ref{lem: T1-1}, 
we have finished the proof of Theorem \ref{thm: mthm2-0} for the case {\bf (A)}.  \qed

Now we prove Theorem \ref{thm: mthm2-0} in the case of\deleted{the} condition {\bf (B)}.
\\
{\it Proof of Theorem \ref{thm: mthm2-0} \added{in the case of {\bf (B)}}}.
By using Markov property, we have that, 
for any Borel measurable functions $F: C([\e,\infty), X) \to \R,$
\begin{align*}
\mathbb E^{\x_n}(F(B^n_{\e+\cdot}))&=\mathbb E^{\x_n}\Bigl(\mathbb E^{B_\e^n}\bigl(F(B^n_{\cdot})\bigr)\Bigr)
\\
&=\int_{X_n}\mathbb E^y\bigl(F(B^n_\cdot)\bigr)p_n({\e},\x_n, dy) 
\\
&=\mathbb E^{p_n({\e},\x_n, dy)}(F(B^n_{\cdot})).
\end{align*}
By Theorem \ref{thm: CONHF}, it holds that $p_n({\e},\x_n, dy) \to p_\infty({\e},\x_\infty, dy)$ in $W_2$-sense and thus also in $W_1$-sense (see e.g., \cite[Remark 6.6]{V09}). Therefore,\deleted{the} condition (i) in Lemma \ref{lem: T1} holds with $\nu_n=p_n({\e},\x_n, dy)$. Moreover,\deleted{the} condition (ii) with $\nu_n=p_n({\e},\x_n, dy)$ in Lemma \ref{lem: T1} also holds by the assumption {\bf (B)}. Therefore, by Lemma \ref{lem: FDC1-2}, and Lemma \ref{lem: T1} with $\nu_n=p_n({\e},\x_n, dy)$, we have the desired result.   \qed
 

 \section{Proof of Theorem \ref{thm: mthm2-1}} \label{sec: proof}
{\it Proof of Theorem \ref{thm: mthm2-1}}:  
The goal of the proof is to show \added{the pmG convergence of $\mathcal X_n$ to $\mathcal X_\infty$,} that \added{is}, for any $f \in C_{bs}(X)$ (recall $C_{bs}(X)$ means the set of bounded continuous functions with bounded supports), we have 
\begin{align*} 
\int_{X}fdm_n \to \int_{X}fdm_\infty \quad \text{as $n \to \infty$.} 
\end{align*}
We first consider the case of $K>0$. 
\\
\underline{The case of $K>0$:}
\\
Let $\lambda_n^1$ be the spectral gap of ${\sf Ch}_n$:
\begin{align*} 
\lambda_n^1 :=\inf\{\frac{{\sf Ch}_n(f)}{\|f\|^2_{L^2(m_n)}}: f \in {\rm Lip}(X_n)\setminus \{0\}, \int_{X_n} f dm_n =0\}.
\end{align*}
The following is a well-known fact (easy to obtain by using the spectral resolution):
\begin{align} \label{eq: equiri}
\|P_t^n  - m_n(\cdot)\|_{2 \to 2} \le e^{-\lambda_n^1 t}, \quad \forall t>0,
\end{align}
whereby $\|\cdot\|_{2 \to 2}$ means the operator norm from $L^2(X_n;m_n)$ to $L^2(X_n;m_n)$, and $m_n(f):=\frac{1}{m_n(X_n)}\int_{X_n}f dm_n$.

By \eqref{eq: equiri} and the assumption \eqref{eq: bound of HK}, we have that, for any $t >t_*$ ($t_*$ appeared in the assumption \eqref{eq: bound of HK}),
\begin{align} \label{ineq: speed of equi}
\|p_n(t,\x_n,\cdot)-\frac{\1}{m_n(X_n)}\|_{L^2(m_n)} &= \|(P_s^n-m_n(\cdot))p_n(t-s,\x_n,\cdot)\|_{L^2(m_n)} \notag
\\
&\le e^{-\lambda_n^1s}\|p_n(t-s, \x_n, \cdot)\|_{L^2(m_n)}  \notag
\\
&=e^{-\lambda_n^1 s}\bigl(p_n(2(t-s), \x_n, \x_n)\bigr)^{1/2}\quad (0 < s< t, \quad t_* <h:=t-s) \notag
\\
&< M^{1/2}e^{-\lambda_n^1(t-h)}.
\end{align}
Since the global Poincar\'e inequality holds under the CD$(K,\infty)$ condition with a positive $K>0$ (see e.g., \cite[Theorem 30.25]{V09}), we have that there exists a positive constant $C_P=C_P(K)$ depending only on $K$ so that  
\begin{align} \label{ineq: POINCARE}
\inf_{n \in \N}\lambda_n>C_P>0.
\end{align}
By the condition of $K>0$, there exists a positive constant $C$ so that $\sup_{n \in \N}m_n(X_n)<C$ (see \cite[Theorem 4.26]{Sturm06}).
Thus, by statement (iii)$_{\ge \e}$, \eqref{ineq: speed of equi} and \eqref{ineq: POINCARE},  we have that, for any $\delta>0$,
\begin{align*} 
&\Bigl| \int_{X}fdm_n -\int_{X}fdm_\infty\Bigr|  \notag
\\
&=\Bigl| \int_{X}fdm_n -m_n(X_n)\mathbb E_n^{\x_n}(f(B_t^n))+m_n(X_n)\mathbb E_n^{\x_n}(f(B_t^n))-m_\infty(X_\infty)\mathbb E_\infty^{\x_\infty}(f(B_t^\infty))
\\
&\quad  \quad +m_\infty(X_\infty)\mathbb E_\infty^{\x_\infty}(f(B_t^\infty))-\int_{X}fdm_\infty\Bigr| \notag
\\
& \le  C\Bigl( \int_{X}|p_n(t,\x_n,y)-\frac{\1}{m_n(X_n)}|fdm_n +|\mathbb E_n^{\x_n}(f(B_t^n))-\mathbb E_\infty^{\x_\infty}(f(B_t^\infty))|
\\
&\quad \quad + \int_{X}|p_\infty(t,\x_n,y)-\frac{\1}{m_\infty(X_\infty)}|fdm_\infty \Bigr) \notag
\\
&\le C\Bigl(\|f\|_{L^2(m_n)}\|p_n(t,\x_n,\cdot)-\frac{\1}{m_n(X_n)}\|_{L^2(m_n)}+|\mathbb E_n^{\x_n}(f(B_t^n)) -\mathbb E_\infty^{\x_\infty}(f(B_t^\infty))|
\\
&\quad \quad +\|f\|_{L^2(m_\infty)}\|p_\infty(t,\x_\infty,\cdot)-\frac{\1}{m_\infty(X_\infty)}\|_{L^2(m_\infty)} \Bigr)\notag
\\
&\le C\Bigl( \|f\|_{L^2(m_n)}Me^{-\lambda_n^1(t-h)}+|\mathbb E_n^{\x_n}(f(B_t^n))-\mathbb E_\infty^{\x_\infty}(f(B_t^\infty))|+\|f\|_{L^2(m_\infty)}Me^{-\lambda_\infty^1(t-h)} \Bigr)
\\
&\le C\Bigl( \|f\|_{L^2(m_n)}Me^{-C_P(t-h)}+|\mathbb E_n^{\x_n}(f(B_t^n))-\mathbb E_\infty^{\x_\infty}(f(B_t^\infty))|+\|f\|_{L^2(m_\infty)}Me^{-C_P(t-h)}\Bigr)
\\
&\to \delta +0 +\delta \quad \text{as $n \to \infty$ and sufficiently large $t$}.
\end{align*}
Thus we finish the proof of Theorem \ref{thm: mthm2-1} for the case of $K>0$.
\\
\underline{The case of $\sup_{n \in \N}{\rm diam}(X_n)<D$:}
\\
The case of $\sup_{n \in \N}{\rm diam}(X_n)<D$ can be proved in the same way as the case of $K>0$ since the local Poincar\'e inequality holds for any RCD$(K,\infty)$ spaces (see \cite[Theorem 1.1]{Raj12b}). If $\sup_{n \in \N}{\rm diam}(X_n)<D$ holds, then the local Poincar\'e inequality \replaced{is equivalent to}{means} the global Poincar\'e inequality and the proof will be the same as the case of $K>0$.
Thus we finish the proof of Theorem \ref{thm: mthm2-1}.
\qed
 
 \section{Proof of Theorem \ref{thm: mthm1}} \label{sec: proof1}
\deleted{In this section, we prove Theorem \ref{thm: mthm1}. In the previous sections, we have already proved the implications (i) $\iff$ (ii) and (iii)$_{\ge 0}$ (or (iii)$_{\ge \e}$) $\implies$ (i) (since under  {Assumption \ref{asmp: 1}}, the conditions assumed in Theorem \ref{thm: mthm2-1} are satisfied (especially condition \eqref{eq: bound of HK} is satisfied by the uniform Gaussian upper heat kernel estimate \eqref{GE}). Thus we only have to show the implication (i) $\implies$ (iii)$_{\ge 0}$.}
\added{In this section, we prove Theorem \ref{thm: mthm1}. In the previous sections, we have already proved (i) $\iff$ (ii) by Theorem \ref{thm: mthm2}. 
If $\sup_{n}{\rm diam}(X_n)<\infty$, then the implication of (iii)$_{\ge \e}$ (or (iii)$_{\ge 0}$) $\implies$ (i) follows from Theorem \ref{thm: mthm2-1}.}
 \added{Thus we only have to show the implication (i) $\implies$ (iii)$_{ \ge 0}$.}
 
 {\it Proof of (i) $\implies$ (iii)$_{ \ge 0}$ in Theorem \ref{thm: mthm1}}.
 Since we have already shown the weak convergence of the laws of finite-dimensional distributions in Lemma \ref{lem: FDC1} for the general RCD$(K,\infty)$ case, what we should prove is only the tightness of the Brownian motions on $C([0,\infty];X)$.
\begin{lem} \label{lem: T}
 $\{\B_n\}_{n \in \N}$ is tight in $\mathcal P(C([0,\infty), X))$.
 \end{lem}
 \proof Since $x_n$ converges to $x_\infty$ in $(X,d)$, \added{the set of }the laws of the initial distributions $\{B^n_0\}_{n \in \N}=\{\delta_{x_n}\}_{n \in \N}$ is clearly tight in $\mathcal P(X)$. 
Thus it suffices to show the following (see \cite[Theorem 12.3]{B99}): 
for each $T>0$, there exist $\beta>0$, $C>0$ and $\theta>1$ such that, for all $n \in \N$
	\begin{align} \label{criterion: tight5}
	&\mathbb E^{x_n}[\tilde{d}^{\beta}(B^{n}_t, B^{n}_{t+h})]
	\le Ch^{\theta}, \quad (0 \le t \le T\ \text{and}\ 0 \le h \le 1), 
	\end{align}
	whereby $\tilde{d}(x,y):=d(x,y) \wedge 1$.

We first give a heat kernel estimate. Let $\sup_{n \in \N}{\rm diam}(X_n)<D.$ By the generalized Bishop--Gromov inequality \cite[{Proposition 3.6}]{EKS15}, we have the following volume growth estimate: \added{for any $D>0$,} there exist positive constants $\nu=\nu(N,K, {D})>0$ and $c=c(N,K,{D})>0$ such that, for all $n \in \N$
\begin{align} \label{VE}
	m_n(B_r(x)) \ge cr^{2\nu} \quad  (0 \le  r \le {D}).
\end{align} 
In fact, by the generalized Bishop--Gromov inequality, we have \added{that, for any $0<r\le D<\infty$,}
\begin{align*}
m_n(B_r(x))\ge \frac{\int_0^r\Theta_{K/N}(t)^Ndt}{\int_0^{{D}}\Theta_{K/N}(t)^Ndt}m_n(B_{{D}}(x)) {\ge}{c(N,K,{D})}\int_0^r\Theta_{K/N}(t)^Ndt.
\end{align*}
Here $c(N,K,{D})=\frac{1}{\int_0^{{D}}\Theta_{K/N}(t)^Ndt}$. Thus we have \eqref{VE}.
Combining \added{\eqref{GE1}} with \eqref{VE}, we have the following uniform upper heat kernel estimate: \deleted{({\bf uniform Gaussian estimate})}
\begin{align} \label{GE}
p_{n}(t,x,y) \le \frac{C_1}{ct^\nu}\exp\Bigl\{-C_2\frac{d_{n}(x,y)^2}{t}\Bigr\},
\end{align}
for all $x, y \in X_{{n}}$ and $0<t\le D^{2}$. Here constants $C_1, C_2, c, \nu$ only depend\deleted{s} on the given constants $N,K,D$. Note that the constant $C_3$ in \eqref{GE1} can be taken as zero under $\sup_{n \in \N}{\rm diam}(X_n)<D$ according to \cite{Sturm96, St98} (note that the MCP condition is satisfied under the assumption of Theorem \ref{thm: mthm1}). 

	Take $\beta>0$ such that $\beta/2-\nu>1$, and set $\theta=\beta/2-\nu$.
By the Markov property, we have 
\begin{align} \label{eq: tight5}
	&\text{L.H.S. of \eqref{criterion: tight5}}  \notag
	\\
	& =\int_{X_n \times X_n} p_n(t,x_n,y)p_n(h,y,z)\tilde{d}^\beta(\iota_n(y),\iota_n(z)) m_n(dy)m_n(dz) \notag
	\\
	& \le \int_{X_n \times X_n} p_n(t,x_n,y)p_n(h,y,z){d}^\beta(\iota_n(y),\iota_n(z))m_n(dy)m_n(dz).
\end{align}
By the Gaussian heat kernel estimate \eqref{GE}, we have
\begin{align} 
	\int_{X_n} &p_n(s,y,z){d}^\beta(\iota_n(y),\iota_n(z)) m_n(dz) \notag
	\\
	&\le \frac{C_1}{cs^{\nu}}\int_{X_n}\exp\Bigl(- C_2\frac{d_n(y,z)^2}{s}\Bigr){d}^\beta(\iota_n(y),\iota_n(z))m_n(dz) \notag
	\\
	&\le \frac{C_1}{cs^{\nu}}\int_{X_n}\exp\Bigl(- C_2\frac{d_n(y,z)^2}{s}\Bigr){d}^\beta_n(y,z)m_n(dz) \notag
	\\
	& \le C_1c^{-1}C_2^{2/\beta}s^{\beta/2-\nu}m_n(X_n)\sup_{y,z\in X_n}\Bigl\{\Bigl(C_2\frac{d_n(y,z)^2}{s}\Bigr)^{\beta/2}\exp\Bigl(- C_2\frac{d_n(y,z)^2}{s}\Bigr)\Bigr\} \notag
	\\ \label{ineq: UGEP}
	& \le  C_1c^{-1}C_2^{2/\beta}C_3M_\beta s^{\beta/2-\nu}\notag
	\\
	& = C_4s^{\beta/2-\nu},
\end{align}
whereby $M_\beta:=\sup_{t \ge 0}{t^{\beta/2}\exp(-t)}$, $C_3=\sup_{n \in \N}m_n(X_n)$ and $C_4=C_4(N,K,D, \beta)=C_1c^{-1}C_2^{2/\beta}C_3M_\beta$ are constants dependent only on $N,K,D, \beta$ (independent of $n$). Note that, since $m_n$ converges weakly to $m_\infty$ and $m_\infty(X_\infty)<\infty$ because of ${\rm diam}(X_\infty)<D$, we have that $\sup_{n \in \N}m_n(X_n)=C_3<\infty$.
By \eqref{ineq: UGEP}, we have
\begin{align*}
	\text{R.H.S. of \eqref{eq: tight5}} & \le C_4h^{\beta/2-\nu} \int_{X_n} p_n(t,x_n,y)m_n(dy) \notag
	\\
	& \le C_4h^{\beta/2-\nu}.
\end{align*}
Thus we finish the proof. 
 \qed

Thus we have completed the proof \added{of} Theorem \ref{thm: mthm1}.
\qed

\section*{Acknowledgment}
The author appreciates Prof.\ Takashi Kumagai for suggesting a constructive comment for the proof in Theorem \ref{thm: mthm2-0} in the case of {\bf (B)}.
He also expresses his great appreciation to Prof.\ Shouhei Honda and Prof.\ Yu Kitabeppu for giving a lot of valuable comments about $\mathrm{RCD}^*(K,N)$ spaces and sparing much time to discuss.
He also appreciates Dr.\ Han Bangxian for fruitful discussions about Theorem \ref{thm: mthm2-1}.
He is indebted to Mrs.\ Anna Katharina \added{Suzuki-}Klasen for her attentive proofreading with helpful suggestions and constructive comments. Finally, he appreciates the anonymous referee for making several corrections. 
This work was supported by Grant-in-Aid for JSPS Fellows Number 261798 and also supported by University Grants for student exchange between universities in partnership of Kyoto University and Hausdorff Center for Mathematics in Bonn.

\end{document}